\documentclass[a4paper,10pt,reqno]{article}
\usepackage{amsmath}
\usepackage{amssymb}
\usepackage{amsthm}
\usepackage[colorlinks=true, bookmarks=true, pdfstartview=FitH, linktocpage=true, linkcolor = magenta, citecolor = blue]{hyperref}
\usepackage{tikz-cd}
\usepackage{multicol}

\hypersetup{
	colorlinks = true,
	linkcolor = magenta,
	citecolor = blue,
}

\newtheorem{theorem}{Theorem}[section]
\newtheorem{lemma}[theorem]{Lemma}
\newtheorem{proposition}[theorem]{Proposition}
\newtheorem{corollary}[theorem]{Corollary}

\newtheorem{remark}[theorem]{Remark}
\newtheorem{remarks}[theorem]{Remarks\/}

\numberwithin{equation}{section}

\DeclareMathOperator{\Div}{div}
\newcommand{\R}{\mathbb R}
\newcommand{\dd}{\,\mathrm{d}}

\begin{document}
	\title{Generalized Caffarelli--Kohn--Nirenberg and related inequalities for Grushin-type gradients}
	\author{
		Trung Hieu Giang
		\footnote{Addresses of Trung Hieu Giang: Department of Mathematical Analysis, Faculty of Mathematics and Physics, Charles University, Praha, Czech Republic;  Institute of Mathematics, Vietnam Academy of Science and Technology, 18 Hoang Quoc Viet, Hanoi, Vietnam. \textit{Email address: trung-hieu.giang@matfyz.cuni.cz}
		} 
		\,
		Nguyen Minh Tri
		\footnote{Addresses of Nguyen Minh Tri: Institute of Mathematics, Vietnam Academy of Science and Technology, 18 Hoang Quoc Viet, Hanoi, Vietnam. \textit{Email address: triminh@math.ac.vn}
		}
		\ and
		Dang Anh Tuan
		\footnote{Addresses of Dang Anh Tuan: University of Sciences, Vietnam National University, 334 Nguyen Trai, Thanh Xuan, Hanoi, Vietnam. \textit{Email address: datuan1105@gmail.com}
		}
		\footnote{Corresponding author.}
	}
	\date{\today}
	\maketitle
	
	\begin{abstract}
		In this paper, we establish necessary and sufficient conditions for a generalized Caffarelli--Kohn--Nirenberg inequality involving the Grushin-type gradient $\nabla_G=(\nabla_x,|x|^\alpha\nabla_y)$, where $\alpha\geq0$, in arbitrary dimensions. When $\alpha=0$, our results recover the classical Caffarelli--Kohn--Nirenberg inequalities and extend the mixed radial and cylindrical weight inequalities of Li and Yan \cite{LY23} to arbitrary dimensions of the $y$ variables. Our approach combines a change of variables with approximation arguments to reduce the Grushin inequalities to Euclidean inequalities with mixed weights. We also characterize the admissible parameters in Hardy--Sobolev inequalities, including the endpoint cases $p=1$ and $q=p$, and derive multiplicative Sobolev and Gagliardo--Nirenberg inequalities. Finally, we establish a weighted Moser--Trudinger inequality on bounded domains and deduce exponential integrability for subcritical growth with arbitrary locally integrable cylindrical power weights. By letting $\alpha=0$, our results also yield the corresponding inequalities for the Euclidean gradient.
	\end{abstract}
	\
	
	
	\
	

	\section{Introduction}
	\label{sec1}
	Weighted inequalities have become a central topic in modern analysis because the presence of weights introduces features that cannot be captured by the Euclidean metric alone. They arise naturally in degenerate elliptic problems in the absence of uniform ellipticity and provide essential tools for determining the natural energy space, the admissible growth of nonlinear terms, the compactness required by variational methods, and, frequently, the critical threshold separating existence from nonexistence. Accordingly, a substantial literature has developed around weighted inequalities and their applications (see, for instance, \cite{CDG93,Tri98pp2,KL12,CCL20} and the references therein).
	
	Let $M,N\geq1$ be integers, let $\alpha \geq 0$, and
	$z=(x,y)\in\R^M\times\R^N$. In this paper, we focus on the Grushin-type gradient
	\begin{equation}
		\label{sec1-eqn1}
		\nabla_G u=(\nabla_xu,|x|^\alpha\nabla_yu)
	\end{equation}
	and the associated homogeneous dimension
	\begin{equation}
		\label{sec1-eqn2}
		Q=M+(1+\alpha)N.
	\end{equation}
	The quadratic operator corresponding to \eqref{sec1-eqn1} is the classical Grushin-type operator $-\Delta_x-|x|^{2\alpha}\Delta_y$. The geometry is generated by the following anisotropic dilations
	\begin{equation}
		\label{sec1-eqn3}
		\delta_\lambda(x,y)=(\lambda x,\lambda^{1+\alpha}y),
		\qquad \lambda>0.
	\end{equation}
	It can be seen that the dilations $\delta_\lambda(x,y)=(\lambda x,\lambda^{1+\alpha} y)$ satisfy $r_\alpha(\delta_\lambda z)=\lambda r_\alpha(z)$ and have Jacobian $\lambda^Q$, where 
	\begin{equation}
		\label{sec1-eqn4}
		r_\alpha(x,y)=\left(|x|^{2(1+\alpha)}+(1+\alpha)^2|y|^2\right)^\frac{1}{2(1+\alpha)}.
	\end{equation}
	
	The Grushin operator is a basic prototype of degenerate elliptic operators. When $\alpha > 0$, it is strongly elliptic away from $\{x=0\}$ and degenerates precisely on this set. Its degeneracy has motivated several research directions, and thus, we only give a brief review here. In particular, the hyperbolic symmetry and spectral estimates were respectively investigated in \cite{Bec00} and \cite{CCL20}. Another important geometric aspect of this structure is the associated isoperimetric problem. Monti and Morbidelli \cite{MM04} proved a sharp isoperimetric inequality in the Grushin plane and determined the corresponding isoperimetric sets. Franceschi and Monti \cite{FM16} subsequently studied the isoperimetric problem in higher-dimensional Grushin spaces and $H$-type groups. The best constants of Sobolev-type embeddings have also been extensively studied in the Grushin setting. Beckner obtained sharp Sobolev estimates in certain low-dimensional cases, while Dou, Sun, Wang, and Zhu established sharp Baouendi--Grushin inequalities and classified the corresponding extremals in a broad parameter range \cite{Bec00,DSWZ22}. In the quadratic case, Monti \cite{Mon06} investigated the existence, symmetry, and uniqueness properties of extremal functions for the critical Sobolev inequality associated with the Grushin gradient. Weighted embeddings involving the Grushin gradient were also established by Nga, Tri, and Tuan \cite{NTT24} in the Grushin plane and by Giang, Tri, and Tuan \cite{GTT26} in a three-dimensional Grushin setting. The results in \cite{NTT24, GTT26} are then applied to prove compact weighted embeddings (see also \cite{LTT24}) and, together with a Pohozaev identity, one can establish the existence and nonexistence theory for semilinear degenerate equations. It is also worth mentioning that the analysis of several equations and boundary value problems associated with certain Grushin-type operators was previously studied in \cite{Tri98pp1, Tri98pp2, TT12}. In the context of general vector fields satisfying H\"ormander's condition, Sobolev embeddings, Harnack inequalities, and regularity theory for nonlinear subelliptic equations were
	developed in \cite{CDG93}.
	
	In this paper, for $1\leq p<\infty$, we denote
	\begin{equation}
		\label{sec1-eqn5}
		p_d^*=
		\begin{cases}
			\dfrac{pd}{d-p},&p<d,\\[1ex]
			+\infty,&p\geq d,
		\end{cases}
	\end{equation}
	where
	$$
	d := M+N.
	$$
	
	One of our main results is a characterization of the admissible parameters in the generalized Caffarelli--Kohn--Nirenberg (CKN) inequality \eqref{sec1-eqn12}. We first establish the corresponding Hardy--Sobolev inequality to cover the endpoint $\theta=1$. The change of variables $\xi=|x|^\alpha x/(1+\alpha)$ reduces the problem to Euclidean inequalities with mixed radial and cylindrical weights. Since its inverse is singular at $\xi=0$ when $\alpha>0$, this reduction requires approximation in the relevant weighted norms. We combine this reduction with the results of Cora, Musina, and Nazarov \cite{CMN26,CMN26pp2} and Li and Yan \cite{LY23}, and give separate arguments for the Hardy--Sobolev endpoints $p=1$ and $q=p$. Our first main result is the following characterization.
	
	\begin{theorem}[Weighted Hardy--Sobolev inequality]
		\label{sec1-thm1}
		Let $0<q<\infty$, $1\leq p<\infty$, and $\alpha_1,\beta_1,\alpha_2,\beta_2\in\R$. Assume that $\beta_2 > -M$. There exists a positive constant $C=C(M,N,\alpha,p,q,\alpha_1,\beta_1,\alpha_2,\beta_2)$ such that the following inequality holds
		\begin{equation}
			\label{sec1-eqn6}
			\left(\int_{\R^{d}}r_\alpha^{\alpha_1}|x|^{\beta_1}|u|^q\dd x\dd y\right)^{1/q} \leq C
			\left(\int_{\R^{d}}r_\alpha^{\alpha_2}|x|^{\beta_2} |\nabla_G u|^p\dd x\dd y\right)^{1/p}
		\end{equation}
		for every $u\in C_c^1(\R^{d})$ if and only if
		\begin{align}
			\label{sec1-eqn7}
			&\beta_1>-M, \\
			\label{sec1-eqn8}
			&\frac{Q+\alpha_1+\beta_1}{q} =\frac{Q+\alpha_2+\beta_2-p}{p}>0,\\
			\label{sec1-eqn9}
			&p\leq q\leq p_d^*,\\
			\label{sec1-eqn10}
			&\frac{\alpha_1}{q}\leq\frac{\alpha_2}{p}.
		\end{align}
		In particular, the homogeneity condition \eqref{sec1-eqn8} determines
		\begin{equation}
			\label{sec1-eqn11}
			q=\frac{p(Q+\alpha_1+\beta_1)}{Q+\alpha_2+\beta_2-p}.
		\end{equation}
	\end{theorem}
	
	We also establish related multiplicative Sobolev and interpolation inequalities. For $M\leq p<Q$, Corollary \ref{sec2-cor2} follows from Theorem \ref{sec1-thm1} by independent dilations in the $y$ variables. It gives an inequality with the weight $|x|^{p\alpha}$ on the left-hand side, separating the norms of $|x|^\alpha\partial_{y_j}u$ while retaining $\nabla_xu$ as one factor. The unweighted inequality in Corollary \ref{sec2-cor3} separates all derivatives $\partial_{x_i}u$ and $|x|^\alpha\partial_{y_j}u$ and yields a corresponding Gagliardo--Nirenberg inequality for $1\leq p<Q$. Its proof uses Theorem \ref{sec1-thm1} for the initial case $M=1$, followed by induction on $M$ at $p=1$ and an application to a power of $|u|$ for $p>1$.
	
	We next consider a generalized version of the celebrated Caffarelli--Kohn--Nirenberg inequality \cite{CKN84}:
	\begin{equation}
		\label{sec1-eqn12}
		\aligned
		\left(\int_{\R^{d}}r_\alpha^{\alpha_1}|x|^{\beta_1}|u|^q\dd x \dd y\right)^{1/q}
		\leq{}& C\left(\int_{\R^{d}}r_\alpha^{\alpha_2}|x|^{\beta_2}|\nabla_G u|^p\dd x\dd y\right)^{\theta/p}\\
		&\times\left(\int_{\R^{d}}r_\alpha^{\alpha_3}|x|^{\beta_3}|u|^s\dd x \dd y\right)^{(1-\theta)/s}.
		\endaligned
	\end{equation}
	For the generalized CKN inequality, let $0<q<\infty$ and $0\leq\theta\leq1$. When $\theta>0$, let $1\leq p<\infty$, and when $\theta<1$, let $0<s<\infty$. At $\theta=0$, the factor involving the Grushin gradient is omitted, and no assumptions on $p,\alpha_2,\beta_2$ are imposed. At $\theta=1$, the second factor on the right-hand side of \eqref{sec1-eqn12} is omitted, and no assumptions on $s,\alpha_3,\beta_3$ are imposed. We assume
	\begin{equation}
		\label{sec1-eqn13}
		\beta_j>-M,\qquad Q+\alpha_j+\beta_j>0,
	\end{equation}
	for $j=1$, also for $j=2$ when $\theta>0$, and for $j=3$ when $\theta<1$. These assumptions ensure local integrability of the weighted measures.
	
	It should be noted that Song and Li \cite[Theorem~1.5]{SL19} established necessary and sufficient conditions for a family of Grushin CKN inequalities with $1<p<Q$. Our next result characterizes the admissible parameters for $0\leq\theta\leq1$, extending their result to a broader family of mixed radial and cylindrical weights and a wider range of exponents. The case $\theta=0$ follows from comparison of the two weights, while $\theta=1$ recovers Theorem \ref{sec1-thm1}. For $0<\theta<1$, we also extend the Euclidean mixed-weight characterization of Li and Yan \cite{LY23} from $y\in\R$ to $y\in\R^N$ for arbitrary $N\geq1$. The proof combines a local cylindrical estimate with a decomposition over dyadic annuli and the parameter reduction of \cite[Lemma~5.3]{LY23}. A change of variables and weighted approximation then transfer this characterization to the Grushin setting.
	
	\begin{theorem}[Necessary and sufficient conditions for the generalized CKN inequality]
		\label{sec1-thm2}
		Let $M,N\geq1$ and assume the exponent and integrability conditions stated above, including both cases $\theta=0$ and $\theta=1$. Then \eqref{sec1-eqn12} holds for every $u\in C_c^1(\R^d)$ if and only if
		\begin{align}
			\label{sec1-eqn14}
			\frac{Q+\alpha_1+\beta_1}{q}
			&=\theta\frac{Q+\alpha_2+\beta_2-p}{p}+(1-\theta)\frac{Q+\alpha_3+\beta_3}{s},\\
			\label{sec1-eqn15}
			\frac{\alpha_1}{q}&\leq\theta\frac{\alpha_2}{p}+(1-\theta)\frac{\alpha_3}{s},
		\end{align}
		and
		\begin{align}
			\label{sec1-eqn20}
			&\frac1q\geq\theta\left(\frac1p-\frac1d\right)+\frac{1-\theta}{s},\\
			\label{sec1-eqn21}
			&\frac{M+\beta_1}{q}\geq\theta\frac{M+\beta_2-p}{p}+(1-\theta)\frac{M+\beta_3}{s}
		\end{align}
		are satisfied, together with the following restriction:
		\begin{equation}
			\label{sec1-eqn22}
			\frac1q\leq\frac{\theta}{p}+\frac{1-\theta}{s}
		\end{equation}
		whenever $\theta=0$ or $\theta=1$, or, for $0<\theta<1$, whenever equality holds in \eqref{sec1-eqn21} or
		\begin{equation}
			\label{sec1-eqn23}
			\frac{Q+\alpha_2+\beta_2-p}{p}=\frac{Q+\alpha_3+\beta_3}{s}.
		\end{equation}
	\end{theorem}
	
	\begin{remarks}
		\label{sec1-rem-ckn}
		\begin{itemize}
			\item[(i)]
			By letting $\alpha=0$ and $\beta_1=\beta_2=\beta_3=0$, Theorem \ref{sec1-thm2} recovers the classical CKN inequality in dimensions $d\geq2$, with $p\geq1$, $q,s>0$, and $0\leq\theta\leq1$, subject to the endpoint conventions above, including the exact parameter conditions \cite{CKN84,LY23}. At $\theta=0$, these conditions reduce to $q=s$ and $\alpha_1=\alpha_3$, and the inequality holds with $C=1$.
			
			\item[(ii)] For $\alpha=0$, Theorem \ref{sec1-thm2} extends the characterization of Li and Yan \cite[Theorem 1.1]{LY23} from $N=1$ to arbitrary $N\geq1$, for $0\leq\theta\leq1$. The case $\alpha>0$ further extends this characterization to Grushin gradients through a change of variables and weighted approximation.
			
			\item[(iii)] For $\alpha=0$, Theorem \ref{sec1-thm1} recovers the mixed-weight Hardy--Sobolev inequality of Cora, Musina, and Nazarov \cite[Theorem~1.1]{CMN26pp2} in the range $p>1$ and $q>p$. The same choice in Corollary \ref{sec2-cor3} gives a multiplicative Gagliardo--Nirenberg inequality for the Euclidean gradient. 
			
			{\item[(iv)] For $\alpha>0$, Nga, Tri, and Tuan \cite{NTT24} and Giang, Tri, and Tuan \cite[Theorem~3]{GTT26} established weighted Sobolev inequalities for $\nabla_G$ in the cases $M=N=1$ and $M=2$, $N=1$, respectively. The corresponding inequalities follow from the sufficiency part of Theorem \ref{sec1-thm1} by taking
				$$
				\alpha_1=\alpha_2=\beta_2=0,\qquad \beta_1=p\alpha.
				$$
				When $M=N=1$, this gives
				$$
				q=\frac{p[(p+1)\alpha+2]}{\alpha+2-p},\qquad 1\leq p<\alpha+2.
				$$
				When $M=2$, $N=1$, and $p=2$, we obtain $q=6$.
				
				\item[(v)] The weighted Hardy--Sobolev inequality of Song and Li \cite[Theorem~1.6]{SL19}, proved for $1<p<Q$, is also a special case of Theorem \ref{sec1-thm1}. Indeed, for $0\leq t\leq p$ and $a>1-Q/p$, take
				$$
				q=\frac{p(Q-t)}{Q-p},\qquad \alpha_1=aq-(1+\alpha)t,\qquad \beta_1=\alpha t,\qquad\alpha_2=ap,\qquad \beta_2=0.
				$$
				Their CKN inequality \cite[Theorem~1.5]{SL19} is likewise contained in Theorem \ref{sec1-thm2}.
				
				\item[(vi)] By letting $\alpha_1 = \alpha_2 = \beta_1 = \beta_2 = 0$, Theorem \ref{sec1-thm1} reduces to the unweighted Sobolev-type inequality for Grushin-type gradients, which has been studied by many authors (see, for example, \cite{Tri14}).}
		\end{itemize}
	\end{remarks}

	Finally, we focus on a weighted Moser--Trudinger inequality. Throughout this paper, let $\Omega \subset \R^d$ be a nonempty bounded smooth domain. For $1\leq p<\infty$, we define
	\begin{equation}
		\label{sec1-eqn24}
		S^p_{1,0}(\Omega):= \textrm{ closure of } C_c^\infty(\Omega) \textrm{ with respect to the }{\|\nabla_G\cdot\|_{L^p(\Omega)}}-\textrm{norm}.
	\end{equation}
	Let
	\begin{equation}
		\label{sec1-eqn25}
		Q'=\frac{Q}{Q-1}.
	\end{equation}
	The following result is obtained from the weighted Moser--Trudinger inequality of Cabr\'e and Ros-Oton \cite[Theorem 1.7]{CR13}, using the change of variables in Proposition \ref{sec2-prop1} and H\"older's inequality.
	
	\begin{theorem}[Weighted Moser--Trudinger inequality]
		\label{sec1-thm4}
		Assume that
		$$
		\beta_0>-M,
		0<\nu<Q'.
		$$
		For every $K>0$, there exists a constant $C=C(M,N,\alpha,\beta_0,\nu,K,\Omega)>0$ such that
		\begin{equation}
			\label{sec1-eqn26}
			\sup\limits_{\substack{u\in S^Q_{1,0}(\Omega)\\
					\|\nabla_G u\|_{L^Q(\Omega)}\leq K}}
			\int_\Omega |x|^{\beta_0}\exp\!\left(|u|^\nu\right)\dd x\dd y \leq C.
		\end{equation}
	\end{theorem}
	To prove Theorem \ref{sec1-thm4}, we first establish the following critical weighted estimate. For a measurable set $E\subset\R^d$ and $\gamma>-M$, define
	$$
	|E|_\gamma:=\int_E |x|^\gamma\dd x\dd y.
	$$
	
	\begin{theorem}
		\label{MoTru1}
		Assume that
		$$
		\alpha_0>-M,\qquad \alpha_0\geq-N\alpha,
		\qquad D=Q+\alpha_0.
		$$
		For every $K>0$, there exist constants $c_1,C>0$, depending only on $M,N,\alpha,\alpha_0,K$, such that, for every $u\in C_c^\infty(\Omega)$ satisfying
		$$
		\left(\int_\Omega |x|^{\alpha_0}|\nabla_G u|^D\dd x\dd y\right)^{1/D}\leq K,
		$$
		we have
		\begin{equation}\label{MoTru12}
			\int_\Omega \exp\!\left(c_1|u|^{D/(D-1)}\right)
			|x|^{\alpha_0+D\alpha}\dd x\dd y
			\leq C|\Omega|_{\alpha_0+D\alpha}.
		\end{equation}
	\end{theorem}
	
	For $\alpha=0$, Theorems \ref{sec1-thm4} and \ref{MoTru1} give Moser--Trudinger estimates for the Euclidean gradient. If also $M=1$, Theorem \ref{MoTru1} recovers the estimate of Cabr\'e and Ros-Oton \cite[Theorem 1.7]{CR13} for the weight $|x_1|^{\alpha_0}$. Taking $\alpha=\alpha_0=0$ gives the unweighted estimate with the classical exponent $d/(d-1)$.

	The organization of our paper is as follows. Section \ref{sec2} gives the proof of Theorem \ref{sec1-thm1} and derives its multiplicative Sobolev and Gagliardo--Nirenberg consequences. In Section \ref{sec3}, we prove the generalized CKN inequality. Section \ref{sec4} is devoted to the proofs of Theorems \ref{MoTru1} and \ref{sec1-thm4}. Finally, in Section \ref{sec5}, we give some concluding remarks.

	\section{Weighted Hardy--Sobolev inequality}
	\label{sec2}
	
	In this section, we present the proof of Theorem \ref{sec1-thm1}. The method of the proof involves a change of variables that allows us to transform the inequality \eqref{sec1-eqn6} into an inequality involving the usual gradient. In particular, we begin with the following.
	
	\begin{proposition}
		\label{sec2-prop1}
		\begin{itemize}
			\item[(i)] Let $1\leq p<\infty$, $0<q<\infty$, and let $\alpha_1,\beta_1,\alpha_2,\beta_2\in\R$. For $u\in C_c^1(\R^{d})$, let
			\begin{equation}
				\label{sec2-eqn1}
				\xi=\Psi_\alpha(x):=\frac{|x|^\alpha x}{1+\alpha},\,
				\zeta=(\xi,y),\, \textrm{and }
				v(\xi,y)=u(\Psi_\alpha^{-1}(\xi),y).
			\end{equation}
			The function $v$ is continuous with compact support and is $C^1$ on $\{\xi\ne0\}$. For all $(x,y)$, we have
			$$
			r_\alpha(x,y)=(1+\alpha)^{1/(1+\alpha)}|\zeta|^{1/(1+\alpha)},
			$$
			and 
			\begin{equation}
				\label{sec2-eqn2}
				\aligned
				&\det D\Psi_\alpha(x)=(1+\alpha)^{1-M}|x|^{M\alpha}, \\
				&\frac{|x|^\alpha}{1+\alpha}|\nabla_\zeta v(\Psi_\alpha(x),y)|
				\leq|\nabla_G u(x,y)| \leq |x|^\alpha|\nabla_\zeta v(\Psi_\alpha(x),y)|,
				\endaligned
			\end{equation}
			for all $x\ne0$.
			\item[(ii)] Let
			\begin{equation}
				\label{sec2-eqn3}
				\begin{split}
					A_1&=\frac{\alpha_1}{1+\alpha},\qquad B_1=\frac{\beta_1-M\alpha}{1+\alpha},\\
					A_2&=\frac{\alpha_2}{1+\alpha},\qquad B_2=\frac{\beta_2+\alpha(p-M)}{1+\alpha},\qquad \kappa_i=M-1+A_i+B_i\quad(i=1,2).
				\end{split}
			\end{equation}
			Then the following identity and inequalities hold:
			\begin{equation}
				\label{sec2-eqn4}
				\int_{\R^d}r_\alpha^{\alpha_1}|x|^{\beta_1}|u|^q\dd x\dd y =(1+\alpha)^{\kappa_1} \int_{\R^d}|\zeta|^{A_1}|\xi|^{B_1}|v|^q\dd\zeta,
			\end{equation}
			and
			\begin{align}
				\label{sec2-eqn5}
				(1+\alpha)^{\kappa_2-p}\int_{\R^d}|\zeta|^{A_2}|\xi|^{B_2}|\nabla v|^p\dd\zeta
				&\leq
				\int_{\R^d}r_\alpha^{\alpha_2}|x|^{\beta_2}|\nabla_G u|^p\dd x\dd y
				\nonumber\\
				&\leq (1+\alpha)^{\kappa_2}\int_{\R^d}|\zeta|^{A_2}|\xi|^{B_2}|\nabla v|^p\dd\zeta.
			\end{align}
		\end{itemize}
	\end{proposition}
	
	\begin{proof}
		(i) A direct differentiation gives
		$$
		D\Psi_\alpha(x)=\frac{|x|^\alpha}{1+\alpha}\left(I+\alpha \frac{x}{|x|}\otimes \frac{x}{|x|}\right).
		$$
		Thus, the Jacobian easily follows. The formula for $r_\alpha$ is a direct consequence of the identity $|\xi|=|x|^{1+\alpha}/(1+\alpha)$. 
		
		Next, for $x \ne 0$, the chain rule gives
		$$
		\nabla_xu=D\Psi_\alpha(x)^T\nabla_\xi v,
		$$
		which immediately implies the inequalities in \eqref{sec2-eqn2}.
		
		\noindent(ii) For $\varepsilon>0$, let
		$$
		U_\varepsilon=\{(x,y):|x|>\varepsilon\} \textrm{ and }V_\varepsilon=\{(\xi,y):|\xi|>\varepsilon^{1+\alpha}/(1+\alpha)\}.
		$$
		The transformation is a $C^1$-diffeomorphism between these domains, and
		\begin{equation}
			\label{sec2-eqn6}
			\dd x\dd y=(1+\alpha)^{M-1}|x|^{-M\alpha}\dd\zeta,\qquad |x|=((1+\alpha)|\xi|)^{1/(1+\alpha)}.
		\end{equation}
		The change of variables on $U_\varepsilon$ gives
		$$
		\int_{U_\varepsilon}r_\alpha^{\alpha_1}|x|^{\beta_1}|u|^q\dd x\dd y =(1+\alpha)^{M-1+(\alpha_1+\beta_1-M\alpha)/(1+\alpha)} \int_{V_\varepsilon}|\zeta|^{A_1}|\xi|^{B_1}|v|^q\dd\zeta.
		$$
		By letting $\varepsilon\to0$ and applying the Monotone Convergence Theorem, \eqref{sec2-eqn4} follows.
		
		Next, raising the pointwise gradient comparison in \eqref{sec2-eqn2} to the power $p$ and multiplying by
		$r_\alpha^{\alpha_2}|x|^{\beta_2}$, which is positive on $U_\varepsilon$, gives
		$$
		\aligned
		(1+\alpha)^{-p}\int_{U_\varepsilon}r_\alpha^{\alpha_2}|x|^{\beta_2+p\alpha}|\nabla v(\Psi_\alpha(x),y)|^p\dd x\dd y
		&\leq\int_{U_\varepsilon}r_\alpha^{\alpha_2}|x|^{\beta_2}|\nabla_G u|^p\dd x\dd y\\
		&\leq \int_{U_\varepsilon}r_\alpha^{\alpha_2}|x|^{\beta_2+p\alpha}|\nabla v(\Psi_\alpha(x),y)|^p\dd x\dd y.
		\endaligned
		$$
		Thanks to \eqref{sec2-eqn6}, we have
		$$
		\aligned
		\int_{U_\varepsilon}r_\alpha^{\alpha_2}|x|^{\beta_2+p\alpha}|\nabla v(\Psi_\alpha(x),y)|^p\dd x\dd y  &=(1+\alpha)^{M-1}\int_{V_\varepsilon} r_\alpha^{\alpha_2}|x|^{\beta_2+\alpha(p-M)}|\nabla v|^p\dd\zeta\\
		&=(1+\alpha)^{\kappa_2}\int_{V_\varepsilon}|\zeta|^{A_2}|\xi|^{B_2}|\nabla v|^p\dd\zeta.
		\endaligned
		$$
		Again, by applying the Monotone Convergence Theorem, we obtain \eqref{sec2-eqn5}. Our proof is complete.
	\end{proof}
	
	For later use, we also record the following local integrability conditions: for $\zeta=(\xi,y)\in\R^M\times\R^N$,
	\begin{equation}
		\label{sec2-eqn7}
		|\zeta|^A|\xi|^B\in L^1_{\mathrm{loc}}(\R^d)
		\textrm{ if and only if } M+B>0 \textrm{ and } d+A+B>0.
	\end{equation}
	For arbitrary $a,b\in\R$, the definition of $r_\alpha$ in \eqref{sec1-eqn4} and the change of variables give
	\begin{equation}
		\label{sec2-eqn8}
		r_\alpha^a|x|^b\in L^1_{\mathrm{loc}}(\R^d) \textrm{ if and only if }
		b>-M \textrm{ and } Q+a+b>0.
	\end{equation}
	
	It should be noted that the function $v$ defined in \eqref{sec2-eqn1} is not, in general, continuously differentiable on all of $\R^d$. Thus, the following lemma will be useful in our analysis.
	
	\begin{lemma}
		\label{sec2-lem1}
		Let $u\in C_c^1(\R^{d})$ and let $v$ be defined as in \eqref{sec2-eqn1}. There exist $v_\varepsilon\in C_c^1(\R^d)$, $0<\varepsilon\leq1$, with supports contained in a common fixed compact set, such that $v_\varepsilon\to v$ pointwise. For every $\ell>0$ and $A,B\in\R$ with $M+B>0$ and $d+A+B>0$, the sequence $\{v_\varepsilon\}$ satisfies
		\begin{equation}
			\label{sec2-eqn9}
			\int_{\R^d}|\zeta|^A|\xi|^B|v_\varepsilon-v|^\ell\dd\zeta\longrightarrow0.
		\end{equation}
		If $1\leq p<\infty$, $\beta_2>-M$, and $Q+\alpha_2+\beta_2>0$, the above sequence also satisfies
		\begin{equation}
			\label{sec2-eqn10}
			\int_{\R^d}|\zeta|^{A_2}|\xi|^{B_2}
			|\nabla v_\varepsilon-\nabla v|^p\dd\zeta\longrightarrow0,
		\end{equation}
		where $A_2$ and $B_2$ are given in \eqref{sec2-eqn3}.
	\end{lemma}
	
	\begin{proof}
		Define
		$$
		T_\varepsilon(\xi) =(1+\alpha)^{\frac{1}{1+\alpha}}(|\xi|^2+\varepsilon^2)^{-\frac{\alpha}{2(1+\alpha)}}\xi \, \textrm{ and } \, v_\varepsilon(\xi,y)=u(T_\varepsilon(\xi),y).
		$$
		Since
		$$
		\Psi_\alpha^{-1}(\xi)=(1+\alpha)^{1/(1+\alpha)}|\xi|^{-\frac{\alpha}{1+\alpha}}\xi \textrm{ for } \xi\ne0,
		$$
		the pointwise convergence holds, including at $\xi=0$. For $|\xi|\geq1$ and $\varepsilon\leq1$, we have
		$$|T_\varepsilon(\xi)|\geq (1+\alpha)^{\frac{1}{1+\alpha}}2^{-\frac{\alpha}{2(1+\alpha)}}|\xi|^{\frac{1}{1+\alpha}}.
		$$
		Thus, the compact support of $u$ gives a common compact support for all $v_\varepsilon$. The convergence in \eqref{sec2-eqn9} is a direct consequence of \eqref{sec2-eqn7} and the Dominated Convergence Theorem.
		
		Next, for $\xi\ne0$,
		$$
		DT_\varepsilon(\xi)=(1+\alpha)^{\frac{1}{1+\alpha}}(|\xi|^2+\varepsilon^2)^{-\frac{\alpha}{2(1+\alpha)}} \left(I-\frac{\alpha}{1+\alpha}\frac{\xi\otimes\xi}{|\xi|^2+\varepsilon^2}\right) 
		$$
		and
		$$
		|DT_\varepsilon(\xi)|\leq C|\xi|^{-\frac{\alpha}{1+\alpha}}.
		$$
		The sequence $\{DT_\varepsilon\}_\varepsilon$ converges to $D\Psi_\alpha^{-1}$ pointwise on $\{ \xi \ne 0\}$ as $\varepsilon\to 0$. The chain rule and boundedness of the derivatives of $u$ give, on the common compact support,
		$$
		|\nabla v_\varepsilon|^p+|\nabla v|^p\leq C(1+|\xi|^{-\frac{p\alpha}{1+\alpha}}).
		$$
		To apply the Dominated Convergence Theorem, we check that this bound remains integrable after multiplication by the weight $|\zeta|^{A_2}|\xi|^{B_2}$. By the local integrability criterion \eqref{sec2-eqn7}, this term is locally integrable since
		$$
		M+B_2-\frac{p\alpha}{1+\alpha}=\frac{M+\beta_2}{1+\alpha}>0 \, \textrm{ and } \,d+A_2+B_2-\frac{p\alpha}{1+\alpha}=\frac{Q+\alpha_2+\beta_2}{1+\alpha}>0.
		$$
		Finally, by applying the Dominated Convergence Theorem, we obtain \eqref{sec2-eqn10}. The proof is complete.
	\end{proof}
	
	Proposition \ref{sec2-prop1} and Lemma \ref{sec2-lem1} allow us to transform \eqref{sec1-eqn6} into an inequality involving the usual gradient. More specifically, we have the following.
	
	\begin{corollary}
		\label{sec2-cor1}
		Let $1\leq p, \, q<\infty$ and $\alpha_1,\beta_1,\alpha_2,\beta_2\in\R$. Assume that
		$$
		\beta_1,\beta_2>-M,\qquad Q+\alpha_1+\beta_1>0,\qquad Q+\alpha_2+\beta_2>0.
		$$
		Then the inequality \eqref{sec1-eqn6} for all $u\in C_c^1(\R^d)$ is equivalent, up to constants depending only on the parameters, to
		\begin{equation}
			\label{sec2-eqn11}
			\left(\int_{\R^d}|\zeta|^{A_1}|\xi|^{B_1}|v|^q\dd\zeta\right)^{1/q} \leq C \left(\int_{\R^d}|\zeta|^{A_2}|\xi|^{B_2}|\nabla v|^p\dd\zeta\right)^{1/p}
		\end{equation}
		for all $v\in C_c^1(\R^d)$.
	\end{corollary}
	
	\begin{proof}
		The result follows immediately from Proposition \ref{sec2-prop1} and Lemma \ref{sec2-lem1}. For this reason, the details are omitted.
	\end{proof}
	
	We now present the necessary estimates. The following first-order weighted Sobolev estimate is given in \cite[Section 2.1.7, Corollary 1]{Maz11}. For this reason, its proof is omitted.
	
	\begin{lemma}
		\label{sec2-lem2}
		Let $t\in\R$ and $q$ satisfy
		$$
		1\leq q\leq\frac{d}{d-1},\qquad M+t>0,\qquad q(d-1+t)>N.
		$$
		Then, for every $w\in C_c^1(\R^d)$,
		\begin{equation}
			\label{sec2-eqn12}
			\left(\int_{\R^d}|\xi|^{q(d-1+t)-d}|w|^q\dd\zeta\right)^{1/q}
			\leq C(M,N,q,t)\int_{\R^d}|\xi|^t|\nabla w|\dd\zeta.
		\end{equation}
	\end{lemma}
	
	\begin{lemma}
		\label{sec2-lem3}
		Let $1\leq p<\infty$ and $s,t\in\R$. Assume that $M+t>0$ and $d+t-s-p>0$. Then, for every $v\in C_c^1(\R^d)$,
		\begin{equation}
			\label{sec2-eqn13}
			\left(\frac{d+t-s-p}{p}\right)^p \int_{\R^d}|\xi|^t|\zeta|^{-s-p}|v|^p\dd\zeta \leq \int_{\R^d}|\xi|^t|\zeta|^{-s}|\nabla v|^p\dd\zeta.
		\end{equation}
		If $M+t-p>0$, then, for every $v\in C_c^1(\R^d)$,
		\begin{equation}
			\label{sec2-eqn14}
			\left(\frac{M+t-p}{p}\right)^p \int_{\R^d}|\xi|^{t-p}|v|^p\dd\zeta \leq \int_{\R^d}|\xi|^t|\nabla_\xi v|^p\dd\zeta \leq \int_{\R^d}|\xi|^t|\nabla v|^p\dd\zeta.
		\end{equation}
	\end{lemma}
	
	\begin{proof}
		For $p>1$, \eqref{sec2-eqn13} follows from \cite[Theorem~1]{CMN26}, with cylindrical dimension $M$, $a=t \ne 0$ and $b=s$. When $t=0$, this is the classical radial Hardy inequality. The first inequality in \eqref{sec2-eqn14} is the classical weighted Hardy inequality in the $\xi$ variables \cite{SSW03}.
		
		It remains to give the proof for the case $p=1$. We define
		$$
		\Phi_\tau(v)=\sqrt{v^2+\tau^2}-\tau,\qquad \tau>0.
		$$
		It is easy to see that $\Phi_\tau(v)\in C_c^1(\R^d)$, $0\leq\Phi_\tau(v)\leq|v|$, and $|\nabla\Phi_\tau(v)|\leq|\nabla v|$. Moreover, $\Phi_\tau(v)\to|v|$ pointwise as $\tau\to0$.
		
		To conclude the proof for \eqref{sec2-eqn13}, let
		$$
		W=|\xi|^t|\zeta|^{-s-1} \,\textrm{ and } \,F=\zeta W.
		$$
		For $\xi\neq0$, a direct calculation gives $\Div F=(d+t-s-1)W$. Choose $R$ such that $\operatorname{supp}v\subset \subset B_R$, where $B_R$ denotes the open Euclidean ball of radius $R$ centered at the origin in $\R^d$, and integrate by parts on
		$$
		\Omega_{\varepsilon,\delta}=B_R\cap\{|\xi|>\varepsilon\}\cap\{|\zeta|>\delta\},\qquad 0<\varepsilon<\delta.
		$$
		For the inner cylindrical boundary, notice that the outward normal is $\nu=-(\xi/|\xi|,0)$, thus $F\cdot\nu=-\varepsilon W$. As a consequence, for fixed $\delta>0$, we have
		\begin{equation}
			\label{sec2-eqn15}
			\left|\int_{\{|\xi|=\varepsilon,\ |\zeta|>\delta\}}\Phi_\tau(v)F\cdot\nu\,\dd S\right|\leq C_{v,\delta,R}\varepsilon^{M+t}\longrightarrow0,
		\end{equation}
		uniformly in $\tau$, since $\Phi_\tau(v)\leq|v|$ and $M+t>0$.
		
		After letting $\varepsilon\to 0$, the remaining inner boundary is the sphere $|\zeta|=\delta$. Its outward normal is $-\zeta/\delta$. Consequently,
		\begin{equation}
			\label{sec2-eqn16}
			\left|\int_{\{|\zeta|=\delta\}}\Phi_\tau(v)F\cdot\nu\,\dd S\right|\leq \|v\|_\infty\delta^{d+t-s-1}\int_{\mathbb S^{d-1}}|\omega_\xi|^t\,\dd\omega\longrightarrow0.
		\end{equation}
		Here $\omega_\xi$ denotes the first $M$ components of $\omega\in\mathbb S^{d-1}$. The surface integral is finite thanks to $M+t>0$, and the limit follows from $d+t-s-1>0$.
		
		Again, from \eqref{sec2-eqn7} and the assumption that $M+t>0$, we deduce that $W$ and $|\xi|^t|\zeta|^{-s}$ are locally integrable. It therefore follows from the boundary limits and the Dominated Convergence Theorem that
		$$
		(d+t-s-1)\int_{\R^d}W\Phi_\tau(v)\,\dd\zeta =-\int_{\R^d}F\cdot\nabla\Phi_\tau(v)\,\dd\zeta\leq\int_{\R^d}|\xi|^t|\zeta|^{-s}|\nabla v|\,\dd\zeta.
		$$
		Letting $\tau\to0$ proves \eqref{sec2-eqn13} for $p=1$.
		
		The only thing left is \eqref{sec2-eqn14} when $p=1$. Assume $M+t-1>0$ and let
		$$
		\widetilde{F}=(\xi|\xi|^{t-1},0).
		$$
		By a simple calculation, we obtain
		$$
		\Div \widetilde{F}=(M+t-1)|\xi|^{t-1}\quad\text{for }\xi\neq0.
		$$
		By integrating by parts outside the tube $|\xi|\leq\varepsilon$ and using a similar argument as above, we obtain
		$$
		\aligned
		(M+t-1)\int_{\R^d}|\xi|^{t-1}\Phi_\tau(v)\,\dd\zeta &=-\int_{\R^d}|\xi|^{t-1}\xi\cdot\nabla_\xi\Phi_\tau(v)\,\dd\zeta\\
		&\leq\int_{\R^d}|\xi|^t|\nabla_\xi v|\,\dd\zeta.
		\endaligned
		$$
		Letting $\tau\to0$, we obtain the first inequality in \eqref{sec2-eqn14} at $p=1$. The second follows from $|\nabla_\xi v|\leq|\nabla v|$.
	\end{proof}
	
	\begin{lemma}
		\label{sec2-lem4}
		Let $1\leq p<\infty$, $v\in C_c^1(\R^d)$, and let $s,t\in\R$ satisfy $M+t>0$ and $d+t-s-p>0$. Define $w=|\zeta|^{-s/p}v$ on $\R^d\setminus\{0\}$. There exist functions $w_\varepsilon\in C_c^1(\R^d)$ with common compact support such that $w_\varepsilon\to w$ pointwise on $\R^d\setminus\{0\}$ and
		\begin{equation}
			\label{sec2-eqn17}
			\int_{\R^d}|\xi|^t|\nabla w_\varepsilon-\nabla w|^p\dd\zeta \to 0.
		\end{equation}
		Moreover,
		\begin{equation}
			\label{sec2-eqn18}
			\left(\int_{\R^d}|\xi|^t|\nabla w|^p\dd\zeta\right)^{1/p}\leq\left(1+\frac{|s|}{d+t-s-p}\right)\left(\int_{\R^d}|\xi|^t|\zeta|^{-s}|\nabla v|^p\dd\zeta\right)^{1/p}.
		\end{equation}
	\end{lemma}
	
	\begin{proof}
		For $\zeta\ne0$, a direct computation gives
		$$
		\nabla w=|\zeta|^{-s/p}\nabla v-\frac{s}{p}|\zeta|^{-s/p-2}\zeta v.
		$$
		Then, the triangle inequality in $L^p(\R^d,|\xi|^t\dd\zeta)$ and \eqref{sec2-eqn13} immediately imply \eqref{sec2-eqn18}.
		
		We construct a sequence $\{w_\varepsilon\}$ satisfying \eqref{sec2-eqn17} by a regularizing argument. In particular, choose $\chi\in C^\infty([0,\infty))$ with $\chi=0$ on $[0,1]$, $\chi=1$ on $[2,\infty)$, and $0\leq\chi\leq1$. Let $\chi_\varepsilon(\zeta)=\chi(|\zeta|/\varepsilon)$ and $w_\varepsilon=\chi_\varepsilon w$, extended by zero near the origin. Then, it can be easily verified that $w_\varepsilon\in C_c^1(\R^d)$. The term $(\chi_\varepsilon-1)\nabla w$ converges to zero in the weighted $L^p$ space by the Dominated Convergence Theorem.
		
		For the remaining term, boundedness of $v$ and the assumption $M+t>0$ give
		$$
		\aligned
		\int_{\R^d}|\xi|^t|w\nabla\chi_\varepsilon|^p\dd\zeta&\leq C_v\varepsilon^{-p}\int_{\{\varepsilon<|\zeta|<2\varepsilon\}}|\xi|^t|\zeta|^{-s}\dd\zeta\\
		&\leq C_v\varepsilon^{d+t-s-p}\to0.
		\endaligned
		$$
		This proves \eqref{sec2-eqn17}. Our proof is complete.
	\end{proof}
	
	Now we are in a position to give the proof of Theorem \ref{sec1-thm1}.
	
	\begin{proof}[Proof of Theorem \ref{sec1-thm1}]
		For the reader's convenience, we divide the proof into two steps.
		\vspace{0.2cm}
		
		\noindent \textbf{{Step 1 (Necessity):}} Suppose \eqref{sec1-eqn6} holds for every $C_c^1(\R^d)$ function. Choose a function equal to one near a point $(0,y_0)$ with $y_0\ne0$ and supported away from the origin. The right-hand side integral corresponding to this function is finite thanks to $\beta_2>-M$, so the local integrability of the left-hand side implies $\beta_1>-M$. Next, choose a compactly supported function equal to one on a neighborhood of the origin. Its gradient vanishes there, and the right-hand side integral corresponding to this function is again finite. By \eqref{sec2-eqn8}, finiteness of the left-hand side then gives $Q+\alpha_1+\beta_1>0$.
		
		Let $u$ be a nonzero test function supported away from $\{x=0\}$. For every $\lambda >0$, define $u_\lambda=u\circ\delta_\lambda$. Both sides of \eqref{sec1-eqn6} corresponding to $u$ scale respectively as
		\begin{equation}
			\label{sec2-eqn19}
			\lambda^{-(Q+\alpha_1+\beta_1)/q} \, \textrm{ and } \,\lambda^{-(Q+\alpha_2+\beta_2-p)/p}.
		\end{equation}
		Letting $\lambda$ tend to zero and to infinity proves \eqref{sec1-eqn8}.
		
		To prove $q\geq p$, let $\varphi$ be a nonzero smooth function supported in $\{1<r_\alpha<2\}$ and away from $\{x=0\}$. We define
		$$
		\varphi_j=3^{\frac{j(Q+\alpha_2+\beta_2-p)}{p}}\varphi\circ\delta_{3^j},\qquad j=1,\ldots,n.
		$$
		The supports of the functions $\varphi_j$ lie in pairwise disjoint annuli. Thus, by a simple calculation together with \eqref{sec2-eqn19} and \eqref{sec1-eqn8}, the inequality \eqref{sec1-eqn6} applied to $\sum_{j=1}^n\varphi_j$ gives 
		$$
		C_0n^{1/q}\leq C_1n^{1/p}
		$$ 
		for every $n$, where $C_0$ and $C_1$ are two positive constants independent of $n$. Thus, by letting $n\to \infty$, we obtain $q\geq p$.
		
		Next, let $\mathcal{O}\subset \R^{d}$ be a bounded smooth domain such
		that $\overline{\mathcal O}\cap\{x=0\}=\varnothing$. Then, by restricting the
		inequality \eqref{sec1-eqn6} to the space $C_c^1(\mathcal O)$, it follows that $q\leq p_d^*$. This concludes \eqref{sec1-eqn9}.
		
		Finally, let $e_1$ denote the first coordinate vector in $\mathbb R^N$ and consider the nonzero function $\varphi$ above. Translate it in the $y$-direction by setting $\varphi_R(x,y)=\varphi(x,y-Re_1)$, where $R>0$. On its support, we have
		$$r_\alpha(x,y)\asymp R^{1/(1+\alpha)} \textrm{ as } R\to\infty.$$ Then, it can be easily seen that both sides of \eqref{sec1-eqn6} have respective orders $R^{\alpha_1/((1+\alpha)q)}$ and $R^{\alpha_2/((1+\alpha)p)}$. Their comparison proves \eqref{sec1-eqn10}.
		\vspace{0.2cm}
		
		\noindent \textbf{{Step 2 (Sufficiency):}} Assume \eqref{sec1-eqn7}--\eqref{sec1-eqn10}. By Corollary~\ref{sec2-cor1}, it remains to prove \eqref{sec2-eqn11}. 
		
		To this end, let
		\begin{equation}
			\label{sec2-eqn20}
			s=-\frac{pA_1}{q},\qquad t=B_2+A_2+s,\qquad h=d+A_2+B_2-p.
		\end{equation}
		Our assumptions imply that
		\begin{equation}
			\label{sec2-eqn21}
			M+B_2=\frac{M+\beta_2+\alpha p}{1+\alpha}>0,\qquad
			h=\frac{Q+\alpha_2+\beta_2-p}{1+\alpha}>0,\qquad s+A_2\geq0.
		\end{equation}
		Thus, we have $t\geq B_2$, $M+t>0$, and $d+t-s-p=h$. It also follows from the condition \eqref{sec1-eqn8} that
		\begin{equation}
			\label{sec2-eqn22}
			A_1=-\frac{sq}{p},\qquad B_1=\frac{q(d-p+t)}p-d,\qquad\frac{q(d-p+t)}p>N.
		\end{equation}
		The last inequality is exactly $M+B_1>0$. Since $|\xi|\leq|\zeta|$, we have, for $\xi\ne0$,
		\begin{equation}
			\label{sec2-eqn23}
			|\xi|^t|\zeta|^{-s}=|\xi|^{B_2}|\zeta|^{A_2}\left(\frac{|\xi|}{|\zeta|}\right)^{s+A_2}\leq |\xi|^{B_2}|\zeta|^{A_2}.
		\end{equation}
		
		We first consider the case $p>1$ and $q>p$. We apply
		\cite[Theorem~1.1]{CMN26pp2}, with cylindrical dimension $M$,
		cylindrical weight exponent $a=B_2$, and radial parameters
		$b=-A_2$ and $\gamma=s$. The local integrability assumptions are precisely
		$$M+B_2>0, \, h>0, \textrm{ and } q(d-p+t)/p>N,$$
		which have been verified in \eqref{sec2-eqn21}--\eqref{sec2-eqn22}. Moreover, its two embedding conditions are
		$(d-p)q\leq dp$ and $\gamma\geq b$. The former follows from \eqref{sec1-eqn9} (and is immediate when $p\geq d$), while the latter is $s+A_2\geq0$. Thus, \eqref{sec2-eqn11} follows. 
		
		Next, we turn to the case where $p=1$ and $q>1$. Fix $v\in C_c^1(\R^d)$
		and let $w=|\zeta|^{-s}v$ for $\zeta\ne0$.
		It should be noted that in general, $w$ need not be $C^1$ at the origin.
		However, thanks to Lemma \ref{sec2-lem4}, there exist functions
		$w_\varepsilon\in C_c^1(\R^d)$ converging pointwise to $w$
		away from the origin and whose gradients converge in
		$L^1(\R^d,|\xi|^t\dd\zeta)$.
		
		By \eqref{sec1-eqn9} and \eqref{sec2-eqn22}, we have
		$$
		q\leq d/(d-1) \textrm{ and } q(d-1+t)>N.
		$$
		Thus Lemma~\ref{sec2-lem2} applies to $w_\varepsilon$.
		Fatou's lemma on the left-hand side and convergence of the gradient
		norms on the right-hand side then imply
		$$
		\aligned
		\left(\int_{\R^d}|\zeta|^{A_1}|\xi|^{B_1}|v|^q\dd\zeta\right)^{1/q}
		&=\left(\int_{\R^d}|\xi|^{q(d-1+t)-d}|w|^q\dd\zeta\right)^{1/q}\\
		&\leq C\int_{\R^d}|\xi|^t|\nabla w|\dd\zeta\\
		&\leq C\left(1+\frac{|s|}{h}\right)
		\int_{\R^d}|\xi|^t|\zeta|^{-s}|\nabla v|\dd\zeta\\
		&\leq C\left(1+\frac{|s|}{h}\right)
		\int_{\R^d}|\xi|^{B_2}|\zeta|^{A_2}|\nabla v|\dd\zeta.
		\endaligned
		$$
		Here we used \eqref{sec2-eqn18} and \eqref{sec2-eqn23}.
		This proves \eqref{sec2-eqn11} in this case.
		
		Finally, let $q=p\geq1$. Thanks to \eqref{sec2-eqn22}, we have
		$$
		A_1=-s \textrm{ and } B_1=t-p.
		$$
		Consequently,
		$$
		M+t-p=M+B_1>0.
		$$
		Again, for $w=|\zeta|^{-s/p}v$, consider the approximation sequence $\{w_\varepsilon\}$ from Lemma \ref{sec2-lem4}. By applying the cylindrical Hardy
		inequality \eqref{sec2-eqn14}, passing to the limit by Fatou's lemma and the weighted gradient convergence \eqref{sec2-eqn17}, we obtain
		$$
		\aligned
		\left(\int_{\R^d}|\zeta|^{A_1}|\xi|^{B_1}|v|^p\dd\zeta\right)^{1/p}
		&=\left(\int_{\R^d}|\xi|^{t-p}|w|^p\dd\zeta\right)^{1/p}\\
		&\leq\frac{p}{M+t-p}
		\left(\int_{\R^d}|\xi|^t|\nabla w|^p\dd\zeta\right)^{1/p}\\
		&\leq\frac{p}{M+t-p}\left(1+\frac{|s|}{h}\right)
		\left(\int_{\R^d}|\xi|^t|\zeta|^{-s}|\nabla v|^p\dd\zeta\right)^{1/p}\\
		&\leq\frac{p}{M+t-p}\left(1+\frac{|s|}{h}\right)
		\left(\int_{\R^d}|\xi|^{B_2}|\zeta|^{A_2}|\nabla v|^p\dd\zeta\right)^{1/p}.
		\endaligned
		$$
	\end{proof}
	
	Theorem \ref{sec1-thm1} also yields the following multiplicative inequality, in which the derivatives in the $y$ variables occur as separate factors.
	
	\begin{corollary}[Multiplicative weighted inequality]
		\label{sec2-cor2}
		Let $1\leq p<Q$ and $p\geq M$, and set
		\begin{equation}
			\label{sec2-eqn24}
			q_\alpha=\frac{p(Q+p\alpha)}{Q-p}.
		\end{equation}
		There exists a constant $C=C(M,N,\alpha,p)>0$ such that, for every $u\in C_c^1(\R^d)$,
		\begin{equation}
			\label{sec2-eqn25}
			\left(\int_{\R^d}|x|^{p\alpha}|u|^{q_\alpha}\dd x\dd y\right)^{1/q_\alpha}
			\leq C\|\nabla_xu\|_{L^p(\R^d)}^{\frac{M+p\alpha}{Q+p\alpha}}\prod_{j=1}^N\bigl\||x|^\alpha\partial_{y_j}u\bigr\|_{L^p(\R^d)}^{\frac{1+\alpha}{Q+p\alpha}}.
		\end{equation}
		Equivalently, after changing the constant,
		\begin{equation}
			\label{sec2-eqn26}
			\begin{split}
				\int_{\R^d}|x|^{p\alpha}|u|^{q_\alpha}\dd x\dd y
				\leq{}&C\left(\int_{\R^d}|\nabla_xu|^p\dd x\dd y\right)^{\frac{M+p\alpha}{Q-p}}\\
				&\times\prod_{j=1}^N\left(\int_{\R^d}|x|^{p\alpha}|\partial_{y_j}u|^p\dd x\dd y\right)^{\frac{1+\alpha}{Q-p}}.
			\end{split}
		\end{equation}
	\end{corollary}
	
	\begin{proof}
		Apply Theorem \ref{sec1-thm1} with $\alpha_1=\alpha_2=\beta_2=0$, $\beta_1=p\alpha$, and $q=q_\alpha$. The balance condition \eqref{sec1-eqn8} follows from \eqref{sec2-eqn24}, and $q_\alpha>p$. If $p<d$, direct calculation gives
		$$
		q_\alpha\leq\frac{pd}{d-p} \textrm{ if and only if }
		\alpha(M-p)\leq0,
		$$
		which holds because $p\geq M$. If $p\geq d$, the upper restriction is immediately satisfied. The other hypotheses are also immediate. Thus, we deduce that
		\begin{equation}
			\label{sec2-eqn27}
			\left(\int_{\R^d}|x|^{p\alpha}|u|^{q_\alpha}\dd x\dd y\right)^{1/q_\alpha}
			\leq C\left(A+\sum_{j=1}^N B_j\right)^{1/p},
		\end{equation}
		where, for the fixed function $u$, we write
		$$
		A=\int_{\R^d}|\nabla_xu|^p\dd x\dd y,
		\qquad B_j=\int_{\R^d}|x|^{p\alpha}|\partial_{y_j}u|^p\dd x\dd y.
		$$
		If $u\neq0$, all these quantities are positive. Indeed, $A=0$ makes $u$ independent of $x$, hence zero by compact support. Likewise, $B_j=0$ implies $\partial_{y_j}u=0$ on $\{x\neq0\}$ and then everywhere by continuity, again forcing $u=0$ by compact support. We may therefore assume $A,B_j>0$.
		
		For $\lambda_1,\ldots,\lambda_N>0$, define
		$$
		u_\lambda(x,y)=u(x,\lambda_1y_1,\ldots,\lambda_Ny_N)
		$$
		and let
		$$
		\Lambda=\prod_{j=1}^N\lambda_j.
		$$
		Applying \eqref{sec2-eqn27} to $u_\lambda$ and changing variables in each integral gives
		\begin{equation}
			\label{sec2-eqn28}
			\left(\int_{\R^d}|x|^{p\alpha}|u|^{q_\alpha}\dd x\dd y\right)^{1/q_\alpha}
			\leq C\Lambda^{\frac1{q_\alpha}-\frac1p}
			\left(A+\sum_{j=1}^N\lambda_j^p B_j\right)^{1/p}.
		\end{equation}
		Since
		$$
		\frac1p-\frac1{q_\alpha}=\frac{1+\alpha}{Q+p\alpha},
		$$
		choose $\lambda_j=(A/B_j)^{1/p}$. Substitution into \eqref{sec2-eqn28} yields
		$$
		\left(\int_{\R^d}|x|^{p\alpha}|u|^{q_\alpha}\dd x\dd y\right)^{1/q_\alpha}
		\leq C A^{\frac{M+p\alpha}{p(Q+p\alpha)}}
		\prod_{j=1}^N B_j^{\frac{1+\alpha}{p(Q+p\alpha)}},
		$$
		where we used $Q-N(1+\alpha)=M$ and absorbed $(N+1)^{1/p}$ into the constant. This is in fact \eqref{sec2-eqn25}. Raising both sides to the power $q_\alpha$ gives \eqref{sec2-eqn26}.
	\end{proof}

	For comparison, when there is no weight on the left-hand side, the inequality admits a product involving each $x$ derivative separately. The proof uses induction on $M$, followed by an application to a power of $|u|$.
	
	\begin{corollary}[Multiplicative unweighted and interpolation inequalities]
		\label{sec2-cor3}
		Let $1\leq p<Q$ and $q_0=pQ/(Q-p)$. There exists $C=C(M,N,\alpha,p)>0$ such that
		\begin{equation}
			\label{sec2-eqn29}
			\|u\|_{L^{q_0}(\R^d)}\leq C\prod_{i=1}^M\|\partial_{x_i}u\|_{L^p(\R^d)}^{\frac{Q-N}{MQ}}\prod_{j=1}^N\bigl\||x|^\alpha\partial_{y_j}u\bigr\|_{L^p(\R^d)}^{\frac1Q}
		\end{equation}
		for every $u\in C_c^1(\R^d)$. More generally, if $0<s<\infty$, $0\leq\theta\leq1$, and
		$$
		\frac1q=\theta\left(\frac1p-\frac1Q\right)+\frac{1-\theta}{s},
		$$
		then
		\begin{equation}
			\label{sec2-eqn30}
			\aligned
			\left(\int_{\R^d}|u|^q\dd x\dd y\right)^{1/q}\leq{}& C\prod_{i=1}^M
			\left(\int_{\R^d}|\partial_{x_i}u|^p\dd x\dd y\right)^{\frac{\theta(Q-N)}{MpQ}}\\
			&\times\prod_{j=1}^N\left(\int_{\R^d}|x|^{p\alpha}|\partial_{y_j}u|^p\dd x\dd y\right)^{\frac{\theta}{pQ}}\\
			&\times\left(\int_{\R^d}|u|^s\dd x\dd y\right)^{\frac{1-\theta}{s}}.
			\endaligned
		\end{equation}
	\end{corollary}
	
	\begin{proof}
		We first prove \eqref{sec2-eqn29} for the case $p=1$ by induction on $M$, keeping $N$ and $\alpha$ fixed. When $M=1$, Theorem \ref{sec1-thm1} with all weight exponents zero gives
		$$
		\|u\|_{L^{Q/(Q-1)}(\R^d)}\leq C\left(\|\partial_xu\|_{L^1(\R^d)}
		+\sum_{j=1}^N\bigl\||x|^\alpha\partial_{y_j}u\bigr\|_{L^1(\R^d)}\right).
		$$
		Apply this to $u(x,\lambda_1y_1,\ldots,\lambda_Ny_N)$ and choose
		$\lambda_j=\|\partial_xu\|_{L^1(\R^d)}/\||x|^\alpha\partial_{y_j}u\|_{L^1(\R^d)}$.
		The same change of variables as in \eqref{sec2-eqn28}, now with $1-1/q_0=1/Q$, proves the desired product estimate. For a nonzero compactly supported $C^1$ function, all these derivative norms are positive, as in the proof of Corollary \ref{sec2-cor2}.
		
		Suppose $M\geq2$ and the assertion holds with $M-1$ in place of $M$, whose homogeneous dimension is $Q-1$. Fix $i\in\{1,\ldots,M\}$ and write $t=x_i$ and $\widehat x_i=(x_1,\ldots,x_{i-1},x_{i+1},\ldots,x_M)$. For this proof set
		$$
		a_k=\|\partial_{x_k}u\|_{L^1(\R^d)},\qquad
		b_j=\bigl\||x|^\alpha\partial_{y_j}u\bigr\|_{L^1(\R^d)},\qquad
		I=\int_{\R^d}|u|^{Q/(Q-1)}\dd x\dd y.
		$$
		H\"older's inequality on each slice gives
		$$
		I\leq\int_{\R}
		\left(\int_{\R^{d-1}}|u|^{(Q-1)/(Q-2)}\dd\widehat x_i\dd y\right)^{\frac{Q-2}{Q-1}}
		\left(\int_{\R^{d-1}}|u|\dd\widehat x_i\dd y\right)^{\frac1{Q-1}}\dd t.
		$$
		By the fundamental theorem of calculus and compact support,
		$$\int_{\R^{d-1}}|u(\widehat x_i,t,y)|\dd\widehat x_i\dd y\leq a_i \textrm{ for every } t.
		$$ 
		Apply the induction hypothesis to the first factor, using
		$|\widehat x_i|^\alpha\leq|x|^\alpha$, and then apply H\"older's inequality in $t$. Since the exponents of the slice integrals sum to one, this yields
		$$
		I\leq C a_i^{\frac1{Q-1}}
		\prod_{k\ne i}a_k^{\frac{Q-1-N}{(M-1)(Q-1)}}
		\prod_{j=1}^N b_j^{\frac1{Q-1}}.
		$$
		Multiplying these inequalities for $i=1,\ldots,M$ and taking the $M$th root, we obtain
		\begin{equation}
			\label{sec2-eqn31}
			I\leq C\prod_{i=1}^M a_i^{\frac{Q-N}{M(Q-1)}}
			\prod_{j=1}^N b_j^{\frac1{Q-1}}.
		\end{equation}
		This proves \eqref{sec2-eqn29} for $p=1$.
		
		For $1<p<Q$, put $\kappa=p(Q-1)/(Q-p)>1$ and apply the case $p=1$ to $v=|u|^\kappa\in C_c^1$. Since
		$\kappa Q/(Q-1)=q_0$ and $(\kappa-1)p/(p-1)=q_0$, H\"older's inequality gives
		$$
		\|\partial_{x_i}v\|_{L^1(\R^d)}\leq\kappa\|u\|_{L^{q_0}(\R^d)}^{\kappa-1}\|\partial_{x_i}u\|_{L^p(\R^d)},\qquad\bigl\||x|^\alpha\partial_{y_j}v\bigr\|_{L^1(\R^d)}
		\leq\kappa\|u\|_{L^{q_0}(\R^d)}^{\kappa-1}\bigl\||x|^\alpha\partial_{y_j}u\bigr\|_{L^p(\R^d)}.
		$$
		The exponents on the right-hand side of \eqref{sec2-eqn29} sum to one. Cancelling $\|u\|_{L^{q_0}(\R^d)}^{\kappa-1}$ for $u\ne0$ proves the assertion for $p>1$. The case $u=0$ is trivial. 
		
		Next, let $s>0$ and $0<\theta<1$. The identity
		$$
		\frac{\theta q}{q_0}+\frac{(1-\theta)q}{s}=1
		$$
		shows that $q_0/(\theta q)$ and $s/((1-\theta)q)$ are conjugate exponents greater than one. Hence H\"older's inequality, applied to $|u|^q=|u|^{\theta q}|u|^{(1-\theta)q}$, gives
		$$
		\left(\int_{\R^d}|u|^q\dd x\dd y\right)^{1/q}
		\leq
		\left(\int_{\R^d}|u|^{q_0}\dd x\dd y\right)^{\theta/q_0}
		\left(\int_{\R^d}|u|^s\dd x\dd y\right)^{(1-\theta)/s}.
		$$
		Substituting \eqref{sec2-eqn29} proves \eqref{sec2-eqn30}, with each derivative exponent multiplied by $\theta$. The constant can be chosen independently of $s$ and $\theta$, since $C^\theta\leq\max\{1,C\}$. For $\theta=0$, the estimate is the identity with $q=s$ while for $\theta=1$, it is \eqref{sec2-eqn29}. In these cases, factors with zero exponent are omitted.
	\end{proof}
	
	\begin{remark} 
		For $p=M$, we also have the same multiplicative weighted inequality as Corollary \ref{sec2-cor3}. In fact, there exists a constant $C=C(M,N)>0$ such that, for every $u\in C_c^1(\R^d)$,
		\begin{equation}
			\label{sec2-eqn26_1}
			\begin{split}
				\int_{\R^d}|x|^{M\alpha}|u|^{M(M+N)/N}\dd x\dd y
				\leq{}&C\prod_{i=1}^M\left(\int_{\R^d}|\partial_{x_i}u|^M\dd x\dd y\right)^{\frac{1}{N}}\\
				&\times\prod_{j=1}^N\left(\int_{\R^d}|x|^{M\alpha}|\partial_{y_j}u|^M\dd x\dd y\right)^{\frac{1}{N}}.
			\end{split}
		\end{equation}
	\end{remark}

	\section{Generalized Caffarelli--Kohn--Nirenberg inequality}
	\label{sec3}
	
	This section is devoted to the proof of Theorem \ref{sec1-thm2}. We first establish Corollary \ref{sec3-cor1} to reduce the problem to Euclidean mixed weights. Proposition \ref{sec3-prop1} characterizes the admissible parameters following the approach of Li and Yan \cite{LY23}. The main step is a local cylindrical estimate combining the radial inequality of Li and Yan \cite[Lemma 4.1]{LY23} in the $\xi$ variables with Sobolev embedding in the $y$ variables. Then, summation over dyadic annuli and the parameter reduction of \cite[Lemma 5.3]{LY23} prove sufficiency. Finally, the dilation and disjoint support tests prove necessity.
	
	\begin{corollary}
		\label{sec3-cor1}
		Assume the exponent and integrability conditions preceding Theorem \ref{sec1-thm2}, with $0\leq\theta\leq1$. Define $A_1,B_1$ by \eqref{sec2-eqn3}; when $\theta>0$, also define $A_2,B_2$ there. When $\theta<1$, let
		\begin{equation}
			\label{sec3-eqn1}
			A_3=\frac{\alpha_3}{1+\alpha} \textrm{ and } B_3=\frac{\beta_3-M\alpha}{1+\alpha}.
		\end{equation}
		Then \eqref{sec1-eqn12} for all $u\in C_c^1(\R^d)$ is equivalent, up to constants depending only on the parameters, to
		\begin{equation}
			\label{sec3-eqn2}
			\aligned
			\left(\int_{\R^d}|\zeta|^{A_1}|\xi|^{B_1}|v|^q\dd\zeta\right)^{1/q}
			&\leq C\left(\int_{\R^d}|\zeta|^{A_2}|\xi|^{B_2}|\nabla v|^p\dd\zeta\right)^{\theta/p}\\
			&\quad\times\left(\int_{\R^d}|\zeta|^{A_3}|\xi|^{B_3}|v|^s\dd\zeta\right)^{(1-\theta)/s}
			\endaligned
		\end{equation}
		for all $v\in C_c^1(\R^d)$, where $\zeta=(\xi,y)\in\R^M\times\R^N$. At $\theta=0$, the gradient factor is omitted; at $\theta=1$, the lower-order factor is omitted.
	\end{corollary}
	\begin{proof}
		The result again follows immediately from Proposition \ref{sec2-prop1} and Lemma \ref{sec2-lem1}. For this reason, the details are omitted.
	\end{proof}
	
	In order to prove the sufficiency of Theorem \ref{sec1-thm2}, we need the following lemma.
	
	\begin{lemma}
		\label{sec3-lem1}
		Let $M,N\geq1$, $d=M+N$, $1\leq p<\infty$, $q,s>0$, and $0<\theta<1$. Suppose
		\begin{equation}
			\label{sec3-eqn3}
			\frac Mq+b_1>0,\qquad \frac Mp+b_2>0,\qquad \frac Ms+b_3>0,
		\end{equation}
		and
		\begin{equation}
			\label{sec3-eqn4}
			\frac dq+b_1=\theta\left(\frac dp+b_2-1\right)+(1-\theta)\left(\frac ds+b_3\right),\qquad 0\leq\eta:=\frac\theta p+\frac{1-\theta}s-\frac1q\leq\frac\theta d.
		\end{equation}
		Let $U\subset\R^N$ be a bounded connected Lipschitz domain, $0<\rho_1<\rho_2$, and {$D_i=B_{\rho_i}^M\times U$ }($i=1,2$), where $B_{\rho_i}^M$ denotes the open Euclidean ball of radius $\rho_i$ centered at the origin in $\R^M$. Then
		\begin{equation}
			\label{sec3-eqn5}
			\aligned
			\left(\int_{D_1}|\xi|^{qb_1}|v|^{q}\dd\zeta\right)^{1/q}&\leq C \Bigg[\left(\int_{D_2\setminus D_1}|\xi|^{qb_1}|v|^{q}\dd\zeta\right)^{1/q}\\
			&\quad+\left(\int_{D_2}|\xi|^{pb_2}|\nabla v|^{p}\dd\zeta\right)^{\theta/p}
			\left(\int_{D_2}|\xi|^{sb_3}|v|^{s}\dd\zeta\right)^{(1-\theta)/s}\Bigg]
			\endaligned
		\end{equation}
		for every Lipschitz function $v$ on $\overline{D_2}$, where $C>0$ depends only on $M,N,p,q,s,\theta,b_1,b_2,b_3,\rho_1,\rho_2$, and $U$, but is independent of $v$.
	\end{lemma}
	
	\begin{proof} For the reader's convenience, we divide the proof into two steps.
		\vspace{0.2cm}
		
		\noindent \textbf{Step 1:} We first consider $p=1$ and present the proof for $M\geq2$. When $M=1$, the entire argument below is carried out separately on the positive and negative half cylinders, choosing a median on each outer component and applying the radial estimate on the corresponding half line. Combining the estimates for the two half cylinders at the end of this step yields \eqref{sec3-eqn5}. 
		
		Denote $S:=D_2\setminus D_1$. For $M\geq2$, this is a connected Lipschitz domain. Let $m_v$ be a median of $v$ on $S$, that is, a real number satisfying
		$$
		|\{\zeta\in S:v(\zeta)<m_v\}|\leq\frac{|S|}{2} \textrm{ and }
		|\{\zeta\in S:v(\zeta)>m_v\}|\leq\frac{|S|}{2}.
		$$
		Notice that such a number exists, but need not be unique, and the two sets need not have equal measures. Let $w=v-m_v$. We prove that
		\begin{equation}
			\label{sec3-eqn6}
			\aligned
			\int_{S}|w|\dd\zeta\leq C\int_{S}|\nabla v|\dd\zeta,\\
			|m_v|\leq C_t\left(\int_{S}|v|^{t}\dd\zeta\right)^{1/t}, \quad t>0,
			\endaligned
		\end{equation}
		where $C$ depends only on $M, N, \rho_1, \rho_2,$ and $U$ and $C_t$ depends only on $t$ and $S$. To this end, we denote 
		$$
		v_S := \frac{1}{|S|} \int_S v \dd\zeta.
		$$
		The usual Poincar\'e inequality gives
		$$
		\int_S |v - v_S| \dd \zeta \leq C\int_S |\nabla v| \dd \zeta.
		$$
		If $m_v \geq v_S$, on the set $\{v \geq m_v\}$, we have
		$$
		|v-v_S| \geq m_v - v_S.
		$$
		Notice that $|\{v \geq m_v\}| \geq |S|/2$. Hence
		$$
		\frac{|S|}{2} |m_v - v_S| \leq \int_S |v-v_S| \dd\zeta.
		$$
		If $m_v < v_S$, we consider the set $\{v \leq m_v\}$ and obtain the same estimate. Therefore,
		$$
		\aligned
		\int_S |w| \dd\zeta = \int_S |v - m_v|\dd\zeta &\leq \int_S |v-v_S|\dd\zeta + |S||v_S - m_v| \\
		& \leq 3 \int_S |v-v_S|\dd\zeta \leq C \int_S |\nabla v| \dd\zeta.
		\endaligned
		$$
		
		For the second inequality in \eqref{sec3-eqn6}, recall that the definition of $m_v$ gives
		$$
		|\{\zeta\in S:v(\zeta)\geq m_v\}|\geq\frac{|S|}{2},\qquad
		|\{\zeta\in S:v(\zeta)\leq m_v\}|\geq\frac{|S|}{2}.
		$$
		Thus, we have $|v|\geq|m_v|$ on a set of measure at least $|S|/2$. Therefore, for every $t>0$,
		$$
		\frac{|S|}{2}|m_v|^t\leq\int_S|v|^t\dd\zeta.
		$$
		This proves the second inequality in \eqref{sec3-eqn6} with $C_t=(2/|S|)^{1/t}$, which depends only on $t$ and $|S|$. 
		
		Next, let
		$$
		\tau=1-\frac{N\eta}{\theta}.
		$$
		It can be easily verified that
		\begin{equation}
			\label{sec3-eqn7}
			\frac Md\leq\tau\leq1,\qquad
			\frac1q=\theta\tau+\theta(1-\tau)\frac{N-1}{N}+\frac{1-\theta}{s}.
		\end{equation}
		To apply the radial local estimate \cite[Lemma 4.1]{LY23} in dimension $M$, define new parameters $\widehat\theta, \widehat s,$ and $\widehat b$ such that
		$$
		\widehat\theta=\theta\tau,\qquad\frac{1-\widehat\theta}{\widehat s}=\theta(1-\tau)+\frac{1-\theta}{s},\qquad(1-\widehat\theta)\widehat b=\theta(1-\tau)b_2+(1-\theta)b_3.
		$$
		Here \cite[Lemma~4.1]{LY23} is applied with $n=M$, $K=\R^M$, and its parameters $(s,p,q,a,\gamma_1,\gamma_2,\gamma_3)$ replaced by $(q,1,\widehat s,\widehat\theta,b_1,b_2,\widehat b)$, respectively. These parameters satisfy the balance identity
		$$
		\frac Mq+b_1=\widehat\theta(M+b_2-1)
		+(1-\widehat\theta)\left(\frac M{\widehat s}+\widehat b\right).
		$$
		The corresponding integrability follows from \eqref{sec3-eqn3}. The other requirements are equivalent to $d\eta\leq\theta$ and $\eta\geq0$. Hence the assumptions of \cite[Lemma 4.1]{LY23} are satisfied.
		
		For $y\in U$, let
		$$
		H(y)=\int_{B_{\rho_1}^M}|\xi|^{b_2}|w(\xi,y)|\dd\xi.
		$$
		Applying \cite[Lemma 4.1]{LY23} to $w(\cdot,y)$ and using H\"older's inequality give
		\begin{equation}
			\label{sec3-eqn8}
			\aligned	       \left(\int_{B_{\rho_1}^M}|\xi|^{qb_1}|w(\xi,y)|^{q}\dd\xi\right)^{1/q}&\leq C\Bigg[\left(\int_{B_{\rho_2}^M\setminus B_{\rho_1}^M}|\xi|^{qb_1}|w(\xi,y)|^{q}\dd\xi\right)^{1/q}\\
			&\quad+\left(\int_{B_{\rho_2}^M}|\xi|^{b_2}|\nabla_\xi v(\xi,y)|\dd\xi\right)^{\theta\tau}H(y)^{\theta(1-\tau)}\\
			&\qquad\times\left(\int_{B_{\rho_1}^M}|\xi|^{sb_3}|w(\xi,y)|^{s}\dd\xi\right)^{(1-\theta)/s}\Bigg]
			\endaligned
		\end{equation}
		for $\tau < 1$. If $\tau=1$, the factor containing $H$ is omitted.
		
		We next claim that
		\begin{equation}
			\label{sec3-eqn9}
			\int_U\bigl(H+|\nabla_yH|\bigr)\dd y\leq C\int_{D_2}|\xi|^{b_2}|\nabla v|\dd\zeta.
		\end{equation}
		Choose a nonnegative cutoff
		$\chi\in C_c^\infty(B_{\rho_2}^M)$ such that $\chi=1$ on $B_{\rho_1}^M$. Since $M+b_2>0$, integration by parts gives
		$$
		\begin{aligned}
			(M+b_2)\int_{D_2}\chi|\xi|^{b_2}|w|\dd\zeta
			&=-\int_{D_2}\chi|\xi|^{b_2}\xi\cdot\nabla_\xi|w|
			\dd\zeta\\
			&\quad-\int_{D_2}|\xi|^{b_2}|w|\xi\cdot\nabla_\xi\chi
			\dd\zeta.
		\end{aligned}
		$$
		
		Since $\nabla w=\nabla v$ and $\nabla_\xi\chi=0$ on $B_{\rho_1}^M$, we obtain
		$$
		\begin{aligned}
			\int_U H(y)\dd y&=\int_{D_1}|\xi|^{b_2}|w|\dd\zeta\\
			&\leq C\Bigg[\int_{D_2}|\xi|^{b_2+1}|\nabla_\xi v|\dd\zeta+\int_S|\xi|^{b_2}|w|\dd\zeta\Bigg].
		\end{aligned}
		$$
		The first term is bounded by $C\int_{D_2}|\xi|^{b_2}|\nabla v|\dd\zeta$ because $|\xi|\leq\rho_2$. For the second term, the weight $|\xi|^{b_2}$ is bounded above and below by positive constants on $S$. Hence the first inequality in \eqref{sec3-eqn6} yields
		$$
		\int_S|\xi|^{b_2}|w|\dd\zeta\leq C\int_S|w|\dd\zeta\leq C\int_S|\nabla v|\dd\zeta\leq C\int_{D_2}|\xi|^{b_2}|\nabla v|\dd\zeta.
		$$
		Consequently,
		$$
		\int_U H(y)\dd y\leq C\int_{D_2}|\xi|^{b_2}|\nabla v|\dd\zeta.
		$$
		
		Since $|\xi|^{b_2}$ is integrable on $B_{\rho_1}^M$, we deduce that $H\in W^{1,1}(U)$ and
		$$
		|\nabla_yH(y)|\leq\int_{B_{\rho_1}^M}|\xi|^{b_2}|\nabla_yv(\xi,y)|\dd\xi\quad\text{for almost every }y\in U.
		$$
		Integrating this inequality over $U$ and combining it with the preceding estimate proves
		$$
		\int_U\bigl(H+|\nabla_yH|\bigr)\dd y\leq C\int_{D_2}|\xi|^{b_2}|\nabla v|\dd\zeta,
		$$
		which is \eqref{sec3-eqn9}.

		Sobolev embedding on $U$, together with \eqref{sec3-eqn9}, gives, when $N>1$,
		$$
		\left(\int_U H^{N/(N-1)}\dd y\right)^{(N-1)/N}
		\leq C\int_{D_2}|\xi|^{b_2}|\nabla v|\dd\zeta.
		$$
		When $N=1$, the same bound holds for
		$\operatorname*{ess\,sup}_{y\in U}H(y)$.
		
		Raise \eqref{sec3-eqn8} to the power $q$ and integrate over $U$. By \eqref{sec3-eqn7},  H\"older's inequality and the preceding estimate yield
		$$
		\begin{aligned}
			\left(\int_{D_1}|\xi|^{qb_1}|w|^q\dd\zeta\right)^{1/q}
			&\leq C\Bigg[\left(\int_S|\xi|^{qb_1}|w|^q
			\dd\zeta\right)^{1/q}\\
			&\quad+\left(\int_{D_2}|\xi|^{b_2}|\nabla v|
			\dd\zeta\right)^\theta
			\left(\int_{D_1}|\xi|^{sb_3}|w|^s
			\dd\zeta\right)^{(1-\theta)/s}\Bigg].
		\end{aligned}
		$$
		
		Since $|\xi|$ is bounded above and away from zero on $S$, the second inequality in \eqref{sec3-eqn6}, with $t=q$, and the condition $M+qb_1>0$ imply
		$$
		\left(\int_S|\xi|^{qb_1}|w|^q\dd\zeta\right)^{1/q}+|m_v|
		\left(\int_{D_1}|\xi|^{qb_1}\dd\zeta\right)^{1/q}\leq
		C\left(\int_S|\xi|^{qb_1}|v|^q\dd\zeta\right)^{1/q}.
		$$
		Similarly, taking $t=s$ and using $M+sb_3>0$, we obtain
		$$
		\left(\int_{D_2}|\xi|^{sb_3}|w|^s\dd\zeta\right)^{1/s}
		\leq C\left(\int_{D_2}|\xi|^{sb_3}|v|^s
		\dd\zeta\right)^{1/s}.
		$$
		Combining these estimates with $v=w+m_v$ proves \eqref{sec3-eqn5} for $p=1$. 
		\vspace{0.2cm}
		
		\noindent \textbf{Step 2:} Now we turn to the case $p>1$. Let $p'=p/(p-1)$ and define
		\begin{equation}
			\label{sec3-eqn10}
			\kappa=1+q/p' \textrm{ and }
			\bar\theta=\frac{\theta\kappa}{1+\theta(\kappa-1)}.
		\end{equation}
		Since $\kappa>1$, the function $|v|^\kappa$ is Lipschitz on $\overline{D_2}$. We apply the case $p=1$ to this function, with exponents $q/\kappa$, $s/\kappa$, interpolation parameter $\bar\theta$, and weight exponents
		$$
		\kappa b_1,\qquad b_2+(\kappa-1)b_1,\qquad \kappa b_3.
		$$
		
		We first check the hypotheses. It is easy to see that $0<\bar\theta<1$, 
		$$
		\frac{M}{q/\kappa}+\kappa b_1=\kappa\left(\frac Mq+b_1\right)>0,\textrm{ and }\frac{M}{s/\kappa}+\kappa b_3 =\kappa\left(\frac Ms+b_3\right)>0.
		$$
		Also, notice that
		$$
		M+b_2+(\kappa-1)b_1 =\left(\frac Mp+b_2\right) +\frac q{p'}\left(\frac Mq+b_1\right)>0.
		$$
		Moreover,
		$$
		d+b_2+(\kappa-1)b_1-1=\left(\frac dp+b_2-1\right)+(\kappa-1)\left(\frac dq+b_1\right).
		$$
		Together with \eqref{sec3-eqn4} and the definition of $\bar\theta$, this gives the required balance identity
		$$
		\kappa\left(\frac dq+b_1\right)=\bar\theta\bigl(d+b_2+(\kappa-1)b_1-1\bigr)+(1-\bar\theta)\kappa\left(\frac ds+b_3\right).
		$$
		The exponent restriction follows from
		$$
		\bar\theta+\frac{1-\bar\theta}{s/\kappa}-\frac1{q/\kappa}=\frac{\kappa\eta}{1+\theta(\kappa-1)}\in\left[0,\frac{\bar\theta}{d}\right].
		$$
		
		By the chain rule and H\"older's inequality, using $(\kappa-1)p'=q$, we have
		$$
		\int_{D_2}|\xi|^{b_2+(\kappa-1)b_1}|\nabla(|v|^\kappa)|\dd\zeta\leq\kappa\left(\int_{D_2}|\xi|^{qb_1}|v|^q\dd\zeta\right)^{(\kappa-1)/q}\left(\int_{D_2}|\xi|^{pb_2}|\nabla v|^p\dd\zeta\right)^{1/p}.
		$$
		Thus the case $p=1$ gives
		$$
		\begin{aligned}
			\left(\int_{D_1}|\xi|^{qb_1}|v|^q\dd\zeta\right)^{\kappa/q}&\leq C\Bigg[\left(\int_{D_2\setminus D_1}|\xi|^{qb_1}|v|^q\dd\zeta\right)^{\kappa/q}\\
			&\quad+\left(\int_{D_2}|\xi|^{qb_1}|v|^q\dd\zeta\right)^{(\kappa-1)\bar\theta/q}\\
			&\qquad\times\left(\int_{D_2}|\xi|^{pb_2}|\nabla v|^p\dd\zeta\right)^{\bar\theta/p}\left(\int_{D_2}|\xi|^{sb_3}|v|^s\dd\zeta\right)^{\kappa(1-\bar\theta)/s}\Bigg].
		\end{aligned}
		$$
		We decompose
		$$
		\int_{D_2}|\xi|^{qb_1}|v|^q\dd\zeta=\int_{D_1}|\xi|^{qb_1}|v|^q \dd\zeta+\int_{D_2\setminus D_1}|\xi|^{qb_1}|v|^q\dd\zeta.
		$$
		Using $(a+b)^r\leq C_r(a^r+b^r)$ with $r=(\kappa-1)\bar\theta/q$, we obtain
		$$
		\left(\int_{D_2}|\xi|^{qb_1}|v|^q
		\dd\zeta\right)^{(\kappa-1)\bar\theta/q}\leq C\Bigg[\left(\int_{D_1}|\xi|^{qb_1}|v|^q\dd\zeta\right)^{(\kappa-1)\bar\theta/q}+\left(\int_{D_2\setminus D_1}|\xi|^{qb_1}|v|^q\dd\zeta\right)^{(\kappa-1)\bar\theta/q}\Bigg].
		$$
		Notice that $0<(\kappa-1)\bar\theta<\kappa$. Substituting this into the preceding estimate and applying Young's inequality with conjugate exponents
		$$
		\frac{\kappa}{(\kappa-1)\bar\theta} \textrm{ and } \frac{\kappa}{\kappa-(\kappa-1)\bar\theta},
		$$
		we deduce that
		$$
		\begin{aligned}
			\left(\int_{D_1}|\xi|^{qb_1}|v|^q\dd\zeta\right)^{\kappa/q} &\leq \frac12\left(\int_{D_1}|\xi|^{qb_1}|v|^q\dd\zeta\right)^{\kappa/q}+C\Bigg[\left(\int_{D_2\setminus D_1}|\xi|^{qb_1}|v|^q\dd\zeta\right)^{\kappa/q}\\
			&\quad+\left(\int_{D_2}|\xi|^{pb_2}|\nabla v|^p\dd\zeta\right)^{\kappa\theta/p}\left(\int_{D_2}|\xi|^{sb_3}|v|^s
			\dd\zeta\right)^{\kappa(1-\theta)/s}\Bigg],
		\end{aligned}
		$$
		where we used
		$$
		\frac{\bar\theta}{\kappa-(\kappa-1)\bar\theta}=\theta,\qquad
		\frac{\kappa(1-\bar\theta)}
		{\kappa-(\kappa-1)\bar\theta}=1-\theta.
		$$
		Moving the first term on the right to the left and taking the $\kappa$th root gives
		$$
		\aligned
		\left(\int_{D_1}|\xi|^{qb_1}|v|^q\dd\zeta\right)^{1/q}
		&\leq C\Bigg[\left(\int_{D_2\setminus D_1}|\xi|^{qb_1}|v|^q\dd\zeta\right)^{1/q}\\
		&\quad+\left(\int_{D_2}|\xi|^{pb_2}|\nabla v|^p\dd\zeta\right)^{\theta/p}
		\left(\int_{D_2}|\xi|^{sb_3}|v|^s
		\dd\zeta\right)^{(1-\theta)/s}\Bigg].
		\endaligned
		$$
	\end{proof}
	
	We have the following, which is an extension of \cite[Theorem 1.1]{LY23} to the case $y\in\R^N$ for arbitrary $N\geq1$.
	
	\begin{proposition}
		\label{sec3-prop1}
		Let $M,N\geq1$, $d=M+N$, $0<q<\infty$, and $0\leq\theta\leq1$. When $\theta>0$, let $1\leq p<\infty$, and when $\theta<1$, let $0<s<\infty$. At $\theta=0$, no assumptions are needed on $p,a_2,b_2$. At $\theta=1$, no assumptions are needed on $s,a_3,b_3$. Assume
		\begin{equation}
			\label{sec3-eqn11}
			\begin{gathered}
				M/q+b_1>0,\quad M/p+b_2>0,\quad M/s+b_3>0,\\
				d/q+a_1+b_1>0,\quad d/p+a_2+b_2>0,\quad d/s+a_3+b_3>0,
			\end{gathered}
		\end{equation}
		where the conditions involving $p,a_2,b_2$ are required only when $\theta>0$, and those involving $s,a_3,b_3$ only when $\theta<1$. Then
		\begin{equation}
			\label{sec3-eqn12}
			\begin{split}
				\left(\int_{\R^d}|\zeta|^{qa_1}|\xi|^{qb_1}|v|^{q}\dd\zeta\right)^{1/q}
				&\leq C\left(\int_{\R^d}|\zeta|^{pa_2}|\xi|^{pb_2}|\nabla v|^{p}\dd\zeta\right)^{\theta/p}\\
				&\quad\times\left(\int_{\R^d}|\zeta|^{sa_3}|\xi|^{sb_3}|v|^{s}\dd\zeta\right)^{(1-\theta)/s}
			\end{split}
		\end{equation}
		holds for every $v\in C_c^1(\R^d)$ if and only if
		\begin{align}
			\frac dq+a_1+b_1&=\theta\left(\frac dp+a_2+b_2-1\right)
			+(1-\theta)\left(\frac ds+a_3+b_3\right),\label{sec3-eqn13}\\
			a_1&\leq\theta a_2+(1-\theta)a_3,\qquad
			\eta:=\frac\theta p+\frac{1-\theta}s-\frac1q\leq\frac\theta d,\label{sec3-eqn14}\\
			\frac Mq+b_1&\geq\theta\left(\frac Mp+b_2-1\right)
			+(1-\theta)\left(\frac Ms+b_3\right),\label{sec3-eqn15}
		\end{align}
		and $\eta\geq0$ whenever $\theta=0$ or $\theta=1$, or, for $0<\theta<1$, whenever equality holds in \eqref{sec3-eqn15} or
		\begin{equation}
			\label{sec3-eqn16}
			\frac dp+a_2+b_2-1=\frac ds+a_3+b_3.
		\end{equation}
	\end{proposition}
	
	\begin{proof}
		We divide the proof into two steps. We follow \cite[Section 5]{LY23} for the sufficiency, while the necessity follows from the tests in \cite[Lemma 2.1]{LY23}, adapted to $\xi\in\R^M$ and $y\in\R^N$.
		\vspace{0.2cm}
		
		\noindent\textbf{Step 1 (Sufficiency):}
		If $\theta=0$, the stated conditions reduce to
		$$
		q=s,\qquad a_1+b_1=a_3+b_3,\qquad a_1\leq a_3.
		$$
		Since $|\xi|\leq|\zeta|$, we have, almost everywhere,
		$$
		|\zeta|^{a_1}|\xi|^{b_1}
		=|\zeta|^{a_3}|\xi|^{b_3}
		\left(\frac{|\xi|}{|\zeta|}\right)^{a_3-a_1}
		\leq |\zeta|^{a_3}|\xi|^{b_3},
		$$
		which proves \eqref{sec3-eqn12} with $C=1$.
		
		If $\theta=1$, the conditions give
		$$
		p\leq q,\qquad \frac1q\geq\frac1p-\frac1d,
		\qquad \frac dq+a_1+b_1=\frac dp+a_2+b_2-1>0,
		\qquad a_1\leq a_2.
		$$
		These are the conditions of Theorem \ref{sec1-thm1} with Grushin parameter $\alpha=0$ and weight exponents
		$$
		(\alpha_1,\beta_1,\alpha_2,\beta_2)
		=(qa_1,qb_1,pa_2,pb_2).
		$$
		Thus that theorem gives \eqref{sec3-eqn12}. In this case, \eqref{sec3-eqn15} is immediately satisfied, since the balance identity gives
		$$
		\frac Mq+b_1-\left(\frac Mp+b_2-1\right)
		=N\left(\frac1p-\frac1q\right)+a_2-a_1\geq0.
		$$
		
		Henceforth suppose $0<\theta<1$. We first assume that $\eta\geq0$. Fix $0<\varepsilon<1/4$ and let
		$$
		K=\{(\xi,y):|\xi|\geq\varepsilon|y|\}.
		$$
		The CKN inequality on cones \cite[Theorem 4.1]{LY23}, with weight factors $|\zeta|^{a_j+b_j}$, gives \eqref{sec3-eqn12} restricted to $K$. Its hypotheses follow from \eqref{sec3-eqn13} and $0\leq\eta\leq\theta/d$.
		
		To estimate the complement, define
		\begin{equation}
			\label{sec3-eqn17}
			\begin{split}
				D_{1,k}&=\{(\xi,y):|\xi|<2\varepsilon\,2^k,\quad 2^k<|y|<2^{k+1}\},\\
				D_{2,k}&=\{(\xi,y):|\xi|<4\varepsilon\,2^k,\quad 2^k<|y|<2^{k+1}\},\qquad k\in\mathbb Z.
			\end{split}
		\end{equation}
		These inner cylinders cover $\R^d\setminus K$ up to null sets, and $D_{2,k}\setminus D_{1,k}\subset K$. After dilation by $2^{-k}$ they are fixed product domains with $U=\{1<|y|<2\}$. For $N=1$, we only need to consider its two interval components separately.
		
		On these fixed cylinders $|\zeta|$ is bounded above and below. Define
		$$
		\widetilde b_2=b_2+\frac{\theta a_2+(1-\theta)a_3-a_1}{\theta}\geq b_2.
		$$
		Equation \eqref{sec3-eqn13} gives \eqref{sec3-eqn4} for $b_1,\widetilde b_2,b_3$. All local integrability requirements hold, and $|\xi|^{\widetilde b_2}\leq C|\xi|^{b_2}$ on the fixed outer cylinder. Lemma \ref{sec3-lem1}, followed by a dilation, therefore gives
		\begin{equation}
			\label{sec3-eqn18}
			\aligned			\left(\int_{D_{1,k}}|\zeta|^{qa_1}|\xi|^{qb_1}|v|^{q}\dd\zeta\right)^{1/q}&\leq C\Bigg[\left(\int_{D_{2,k}\setminus D_{1,k}}|\zeta|^{qa_1}|\xi|^{qb_1}|v|^{q}\dd\zeta\right)^{1/q}\\
			&\quad+\left(\int_{D_{2,k}}|\zeta|^{pa_2}|\xi|^{pb_2}|\nabla v|^{p}\dd\zeta\right)^{\theta/p}\left(\int_{D_{2,k}}|\zeta|^{sa_3}|\xi|^{sb_3}|v|^{s}\dd\zeta\right)^{(1-\theta)/s}\Bigg].
			\endaligned
		\end{equation}
		The rescaled domains are fixed, so the dependence of the constant in Lemma \ref{sec3-lem1} on $U$ and the two radii does not introduce any dependence on $k$. The balance identity \eqref{sec3-eqn13} makes the scaling factors on both sides agree. Thus the constant in \eqref{sec3-eqn18} is independent of $k$. 
		
		Thanks to $\eta\geq0$, we have
		$$
		\frac{\theta q}{p}+\frac{(1-\theta)q}{s}\geq1.
		$$
		Consequently, for nonnegative sequences $\{X_k\}$ and $\{Y_k\}$, we have
		$$
		\sum_k X_k^{\theta q/p}Y_k^{(1-\theta)q/s}\leq\left(\sum_kX_k\right)^{\theta q/p}\left(\sum_kY_k\right)^{(1-\theta)q/s}.
		$$
		It follows from this inequality that
		$$
		\aligned
		&\sum_k\left(\int_{D_{2,k}}|\zeta|^{pa_2}|\xi|^{pb_2}|\nabla v|^{p}\dd\zeta\right)^{\theta q/p}\left(\int_{D_{2,k}}|\zeta|^{sa_3}|\xi|^{sb_3}|v|^{s}\dd\zeta\right)^{(1-\theta)q/s}\\
		& \leq \left(\sum_k\int_{D_{2,k}}|\zeta|^{pa_2}|\xi|^{pb_2}|\nabla v|^{p}\dd\zeta\right)^{\theta q/p}\left(\sum_k\int_{D_{2,k}}|\zeta|^{sa_3}|\xi|^{sb_3}|v|^{s}\dd\zeta\right)^{(1-\theta)q/s}.
		\endaligned
		$$
		Since $D_{2,k}\setminus D_{1,k}\subset K$ and the sets $D_{2,k}$ are disjoint up to null sets, we deduce that
		$$
		\sum_k\int_{D_{2,k}\setminus D_{1,k}}|\zeta|^{qa_1}|\xi|^{qb_1}|v|^{q}\dd\zeta
		$$
		is bounded by the estimate on $K$. Raising \eqref{sec3-eqn18} to the power $q$, summing over $k$, and using $(a+b)^q\leq C_q(a^q+b^q)$ for $a,b\geq0$ together with the preceding two estimates prove \eqref{sec3-eqn12} in this case.
		
		For $\eta<0$, we follow \cite[Lemma 5.3]{LY23} and its application to \cite[Theorem 1.1]{LY23}. Let
		$$
		E=\frac dp+a_2+b_2-1,\qquad L=\frac ds+a_3+b_3,\qquad
		T=\frac dq+a_1+b_1=\theta E+(1-\theta)L.
		$$
		The assumptions give $E\ne L$ and strict inequality in \eqref{sec3-eqn15}. For $v\ne0$, multiplication by a positive constant and dilation therefore allow us to assume that both right-hand integrals in \eqref{sec3-eqn12} equal one. For $t$ close to $\theta$, choose
		$$
		\begin{aligned}
			&\frac1{q_t}=\frac tp+\frac{1-t}s,\qquad a_t=tE+(1-t)L-\frac d{q_t}-b_t,\\
			&\max\left\{-\frac M{q_t},\,t(b_2-1)+(1-t)b_3\right\}
			<b_t<b_1+M\left(\frac1q-\frac1{q_t}\right).
		\end{aligned}
		$$
		This interval is nonempty by \eqref{sec3-eqn11} and the strict inequality in \eqref{sec3-eqn15}. For $t$ sufficiently close to $\theta$, we have $q_t>q$, and the parameters $q_t,t,a_t,b_t$ satisfy all the hypotheses of the case $\eta=0$. Choose $t_-,t_+$ in this neighborhood so that
		$$
		t_-E+(1-t_-)L<T<t_+E+(1-t_+)L.
		$$
		Applying the case $\eta=0$ with $t=t_\pm$, under the normalization above, gives
		$$
		\left(\int_{\R^d}|\zeta|^{q_{t_\pm}a_{t_\pm}}|\xi|^{q_{t_\pm}b_{t_\pm}} |v|^{q_{t_\pm}}\dd\zeta\right)^{1/q_{t_\pm}}\leq C.
		$$
		We now estimate the integral on the left-hand side of \eqref{sec3-eqn12}. On $\{|\zeta|<1\}$, H\"older's inequality gives
		$$
		\begin{aligned}
			\int_{|\zeta|<1}|\zeta|^{qa_1}|\xi|^{qb_1}|v|^q\dd\zeta\leq C\left(\int_{|\zeta|<1}|\zeta|^{\frac{qq_{t_-}(a_1-a_{t_-})}{q_{t_-}-q}}|\xi|^{\frac{qq_{t_-}(b_1-b_{t_-})}{q_{t_-}-q}}\dd\zeta\right)^{1-q/q_{t_-}}.
		\end{aligned}
		$$
		The integral on the right is finite since
		$$
		b_{t_-}<b_1+M\left(\frac1q-\frac1{q_{t_-}}\right),\qquad t_-E+(1-t_-)L<T.
		$$
		The same argument on $\{|\zeta|\geq1\}$, with $t_+$ in place of $t_-$, gives a finite bound since
		$$
		b_{t_+}<b_1+M\left(\frac1q-\frac1{q_{t_+}}\right),\qquad t_+E+(1-t_+)L>T.
		$$
		Adding these estimates bounds the left-hand side of \eqref{sec3-eqn12} under the normalization above. Undoing the multiplication and dilation, using \eqref{sec3-eqn13}, proves \eqref{sec3-eqn12}.
		\medskip
		
		\noindent\textbf{Step 2 (Necessity):} The necessity follows by adapting the proof of \cite[Lemma 2.1]{LY23} to $\zeta=(\xi,y)\in\R^M\times\R^N$. Dilation, translation in a fixed $y$ direction, and concentration away from $\{\xi=0\}$ give \eqref{sec3-eqn13} and \eqref{sec3-eqn14}. Scaling the $\xi$ variables gives \eqref{sec3-eqn15}, with transverse dimension $M$. Finally, the same disjoint-support arguments yield $\eta\geq0$ when $\theta=0$, $\theta=1$, or, for $0<\theta<1$, when equality holds in \eqref{sec3-eqn15} or \eqref{sec3-eqn16} holds. 
	\end{proof}

	\begin{proof}[Proof of Theorem \ref{sec1-thm2}]
		By Corollary \ref{sec3-cor1}, it suffices to apply Proposition \ref{sec3-prop1} with
		\begin{equation}\label{sec3-eqn21}
			(a_1,b_1)=(A_1/q,B_1/q),\qquad
			(a_2,b_2)=(A_2/p,B_2/p),\qquad
			(a_3,b_3)=(A_3/s,B_3/s),
		\end{equation}
		omitting the second pair at $\theta=0$ and the third at $\theta=1$. For $j=1$ and, when present, $j=3$, we have
		$$
		M+B_j=\frac{M+\beta_j}{1+\alpha},\qquad
		d+A_j+B_j=\frac{Q+\alpha_j+\beta_j}{1+\alpha}.
		$$
		When $\theta>0$,
		\begin{equation}\label{sec3-eqn22}
			M+B_2-p=\frac{M+\beta_2-p}{1+\alpha},\qquad
			d+A_2+B_2-p=\frac{Q+\alpha_2+\beta_2-p}{1+\alpha}.
		\end{equation}
		In particular, \eqref{sec1-eqn13} implies the applicable conditions in \eqref{sec3-eqn11}. These identities transform \eqref{sec3-eqn13} into \eqref{sec1-eqn14}, the two conditions in \eqref{sec3-eqn14} into \eqref{sec1-eqn15} and \eqref{sec1-eqn20}, and \eqref{sec3-eqn15} into \eqref{sec1-eqn21}. Equality in \eqref{sec3-eqn15} is preserved, and, for $0<\theta<1$, \eqref{sec3-eqn16} becomes \eqref{sec1-eqn23}. Finally, $\eta\geq0$ is precisely \eqref{sec1-eqn22}, with the same endpoint and equality cases. Proposition \ref{sec3-prop1} therefore gives both necessity and sufficiency for all $0\leq\theta\leq1$.
	\end{proof}

	\section{Weighted Moser--Trudinger inequality}
	\label{sec4}
	In order to prove Theorem \ref{MoTru1}, we need another version of the Moser--Trudinger inequality by Cabr\'e and Ros-Oton \cite[Theorem 1.7]{CR13} as follows:
	\begin{proposition}
		\label{Md_MT11}
		Let $\Omega\subset\R^d$ be a bounded domain. Let 
		$$\beta\ge 0, D_\beta=d+\beta.$$
		Then, for every $K>0$, there exist constants $c_1,C>0$, depending only on $\beta,M,D_\beta,K$, such that, if $v\in C_c^1(\Omega)$ and $\int_\Omega |\xi|^\beta|\nabla v|^{D_\beta}\dd\xi\dd y\leq K$, we have
		\begin{equation}\label{MoTru121}
			\int_\Omega \exp(c_1|v|^{D_\beta/(D_\beta-1)})|\xi|^{\beta}\dd\xi\dd y\leq C|\Omega|_{\beta},
		\end{equation}
		where $\xi = (\xi_1,...,\xi_M)$.
	\end{proposition}
	\begin{proof}
		Note that, for $\beta\ge 0,$
		\begin{equation}
			\label{power}
			\frac{1}{M}\sum_{j=1}^M |\xi_j|^\beta\le |\xi|^\beta_\infty\le |\xi|^\beta\le M^{\beta/2}|\xi|^\beta_\infty\le M^{\beta/2}\sum_{j=1}^M |\xi_j|^\beta,
		\end{equation}
		where $|\xi|_\infty=\max_{1\leq j\leq M}|\xi_j| $, so the assumption $\int_\Omega |\xi|^\beta|\nabla v|^{D_\beta}\dd\xi\dd y\leq K$ gives
		$$\int_\Omega |\xi_j|^\beta |\nabla v|^{D_\beta}\dd\xi\dd y\le K.$$
		Hence, from \cite[Theorem 1.7]{CR13}, there exist constants $c_1,C>0$, depending only on $D_\beta,K$, such that
		$$\int_\Omega \exp(c_1|v|^{D_\beta/(D_\beta-1)})|\xi_j|^{\beta}\dd\xi\dd y\le C|\Omega|_{\beta, j}$$
		where $|\Omega|_{\beta,j}=\int_\Omega|\xi_j|^\beta\dd\xi\dd y$. Using \eqref{power} again, we obtain \eqref{MoTru121}.
	\end{proof}
	
	\begin{proof}[Proof of Theorem \ref{MoTru1}]
		Using the transformation in Proposition \ref{sec2-prop1}, we have 
		$$v(\xi, y)=u((\alpha+1)^{\frac{1}{\alpha+1}}|\xi|^{-\frac{\alpha}{\alpha+1}}\xi, y), \tilde{\Omega}=\{(|x|^\alpha x/(\alpha+1), y): (x, y)\in\Omega\}$$
		where $\beta=N\alpha+\alpha_0$. Then $D=D_\beta, |\Omega|_{\alpha_0+D\alpha} = C_1|\tilde{\Omega}|_\beta,$
		\begin{align*}
			\int_\Omega \exp(c_1|u|^{D/(D-1)})|x|^{\alpha_0+D\alpha}\dd x\dd y & =C_1\int_{\tilde{\Omega}} \exp(c_1|v|^{D_\beta/(D_\beta-1)})|\xi|^\beta \dd\xi\dd y,\\
			\int_\Omega |x|^{\alpha_0}|\nabla_G u|^{D}\dd x\dd y & \ge C_2\int_{\tilde{\Omega}} |\xi|^\beta |\nabla v|^{D_\beta} \dd\xi\dd y
		\end{align*}
		where $C_1,C_2>0$ are constants depending on $\alpha,\beta,M,D$. By the hypothesis,
		$$
		\int_{\tilde{\Omega}}|\xi|^\beta|\nabla v|^{D_\beta}\dd\xi\dd y\le C_2^{-1}K^D.
		$$
		Using the approximation in Lemma \ref{sec2-lem1}, we apply Proposition \ref{Md_MT11} with $2C_2^{-1}K^D$ in place of $K$ and then pass to the limit by Fatou's lemma. This proves the result.
	\end{proof}

	\begin{proof}[Proof of Theorem~\ref{sec1-thm4}]
		By approximation in $S^Q_{1,0}(\Omega)$ and Fatou's lemma, it suffices to prove the estimate for $u\in C_c^\infty(\Omega)$. For a function with gradient norm at most $K$, the approximating functions have gradient norms at most $2K$ for all sufficiently large indices; we therefore use the smooth estimate with $2K$ when passing to the limit. Let {$c_1,C$ }be the constants in {Theorem }\ref{MoTru1} for {$\alpha_0=0$. For $\|\nabla_G u\|_{L^Q(\Omega)}\leq K$, }we obtain
		\begin{equation}\label{MoTru212}
			\int_\Omega \exp(c_1|u|^{Q'})|x|^{Q\alpha}\dd x\dd y\leq C|\Omega|_{Q\alpha}.
		\end{equation}
		Let $K_1>0$, $1<r<\infty$, and $r'=r/(r-1)$. Since $\nu<Q'$, we have
		\begin{align}
			\int_{|u|>K_1} &\exp(|u|^\nu)|x|^{\beta_0}\dd x\dd y \le \int_\Omega \exp\left(\frac{1}{K_1^{Q'-\nu}}|u|^{Q'}\right)|x|^{\beta_0}\dd x\dd y\notag\\
			& \le (\int_\Omega|x|^{\beta_1}\dd x\dd y)^{1/r'}(\int_\Omega \exp\left(\frac{r}{K_1^{Q'-\nu}}|u|^{Q'}\right)|x|^{Q\alpha}\dd x\dd y)^{1/r}\label{K_PMT}
		\end{align}
		where $\beta_1=(r\beta_0-Q\alpha)/(r-1)$. Since $\lim_{r\to\infty}\beta_1=\beta_0>-M$, there is an $r>1$ such that $\int_\Omega|x|^{\beta_1}\dd x\dd y$ is finite. For this $r$, let $K_1>0$ such that $K_1^{\nu-Q'}=c_1/r$. Therefore, from \eqref{MoTru212} and \eqref{K_PMT}, we obtain
		\begin{equation}\label{K_PMT1}
			\int_{|u|>K_1}\exp(|u|^\nu)|x|^{\beta_0}\dd x\dd y\le |\Omega|_{\beta_1}^{1/r'}(C|\Omega|_{Q\alpha})^{1/r}.
		\end{equation} 
		On the other hand, since $\nu>0$ and $\beta_0>-M$, the same choice of $K_1$ gives
		\begin{equation}\label{K_PMT2}
			\int_{|u|\le K_1}\exp(|u|^\nu)|x|^{\beta_0}\dd x\dd y\le \exp(K_1^\nu)|\Omega|_{\beta_0}<\infty.
		\end{equation} 
		Combining \eqref{K_PMT1} and \eqref{K_PMT2} proves Theorem \ref{sec1-thm4}.
	\end{proof}
	
	\section{Concluding remarks}
	\label{sec5}
	We have established a weighted Hardy--Sobolev inequality, a generalized Caffarelli--Kohn--Nirenberg inequality, and weighted Moser--Trudinger inequalities for a Grushin-type gradient. Under the stated integrability assumptions, Theorems \ref{sec1-thm1} and \ref{sec1-thm2} give the necessary and sufficient conditions on the parameters of the first two inequalities. Theorem \ref{MoTru1} establishes a weighted Moser--Trudinger inequality for functions satisfying
	$$
	\int_\Omega |x|^{\alpha_0}|\nabla_G u|^D\dd x\dd y\leq K^D,
	$$
	where $D=Q+\alpha_0$, $\alpha_0>-M$, and $\alpha_0\geq-N\alpha$. Its case $\alpha_0=0$, combined with H\"older's inequality, yields the weighted exponential integrability in Theorem \ref{sec1-thm4} for every $\beta_0>-M$ and $0<\nu<Q'$.
	
	The change of variables in Proposition \ref{sec2-prop1} is essential to our approach. It transforms the Grushin gradient into a Euclidean gradient with mixed radial and cylindrical weights. Since the transformed functions need not be continuously differentiable across $\{\xi=0\}$, the reduction also requires approximation given by Lemma \ref{sec2-lem1}.
	
	For the Hardy--Sobolev inequality, the Euclidean result of Cora, Musina, and Nazarov \cite{CMN26pp2} applies when $p>1$ and $q>p$. Section \ref{sec2} treats the remaining cases $p=1$ and $q=p$, including $p=q=1$, through first-order weighted Sobolev and Hardy estimates. The multiplicative consequences in Corollaries \ref{sec2-cor2} and \ref{sec2-cor3} further separate the contributions of individual derivatives. Their proofs use independent dilations in the $y$ variables and induction on the number of $x$ variables.
	
	For the generalized CKN inequality, the Euclidean characterization of Li and Yan \cite{LY23} concerns only a one-dimensional cylindrical axis. Proposition \ref{sec3-prop1} extends this characterization to $\mathbb R^M\times\mathbb R^N$ for arbitrary $M,N\geq1$. The local estimate in Lemma \ref{sec3-lem1} combines their radial inequality in the $\xi$ variables with Sobolev embedding in the $y$ variables. Summation over dyadic regions and the parameter reduction of \cite[Lemma 5.3]{LY23} establish sufficiency, while similar tests as in \cite[Proof of Lemma 2.1]{LY23} give necessity. The change of variables then yields the characterization in the Grushin setting. Theorems \ref{sec1-thm1} and \ref{sec1-thm2} also extend the Hardy--Sobolev and CKN inequalities of Song and Li \cite{SL19} to a broader family of mixed radial and cylindrical weights and a wider range of exponents, including $p=1$.

	\section*{Acknowledgements}
	The first author is supported by the grant PRIMUS/24/SCI/020 of Charles University. All the authors are supported by the project "Some topics in the theory of partial differential equations containing degenerate elliptic operators" at the Institute of Mathematics, VAST, Hanoi.
	
	\section*{Data availability}
	No data were generated or analysed in this study.
	
	\section*{Declaration on the use of {artificial intelligence}}
	The authors used ChatGPT (OpenAI) to assist in further developing their initial mathematical results and proof strategies, and in revising the manuscript. The authors subsequently reworked the arguments and independently verified the mathematical statements, proofs, and references. They take full responsibility for the accuracy, originality, and integrity of the final manuscript.
	
	\section*{Conflicts of interest}
	The authors have no conflicts of interest to declare that are relevant to the content of this article.

	\newcommand \auth{} 
	\newcommand \jour{} 
	\newcommand \book{}

\end{document}